\documentclass[aop,seceqn,MSNbibl,citesort,dvips]{arximspdf}
\usepackage{graphics}
\doi{10.1214/09-AOP506}
\volume{38}
\issue{3}
\pubyear{2010}
\firstpage{1106}
\lastpage{1142}

\makeatletter

\newcommand{\vienas}{{1}}

\newtheorem{theorem}{Theorem}
\newtheorem{prop}{Proposition}[section]
\newtheorem{cor}{Corollary}[section]
\newproclaim{defi}{Definition}[section]
\newtheorem{lemma}{Lemma}[section]

\makeatother

\begin{document}
\begin{frontmatter}

\title{Martin boundary of a killed random walk on a~quadrant}
\runtitle{Martin boundary of a killed random walk on a quadrant}

\begin{aug}
\author[A]{\fnms{Irina} \snm{Ignatiouk-Robert}\corref{}\ead[label=e1]{Irina.Ignatiouk@u-cergy.fr}} and
\author[A]{\fnms{Christophe} \snm{Loree}\ead[label=e2]{Christophe.Loree@u-cergy.fr}}
\runauthor{I. Ignatiouk-Robert and C. Loree}
\affiliation{UMR CNRS 8088, Universite de Cergy-Pontoise}
\address[A]{D\'{e}partement de math\'{e}matiques\\
Universit\'{e} de Cergy-Pontoise\\
2, Avenue Adolphe Chauvin\\
95302 Cergy-Pontoise Cedex\\
France\\
\printead{e1}\\
\phantom{E-mail: }\printead*{e2}} 
\end{aug}

\received{\smonth{2} \syear{2009}}
\revised{\smonth{7} \syear{2009}}

%
\begin{abstract}
A complete representation of the Martin boundary of killed random walks
on the quadrant ${\mathbb N}^*\times{\mathbb N}^*$ is obtained. It is
proved that the
corresponding full Martin compactification of the quadrant
${\mathbb N}^*\times{\mathbb N}^*$ is homeomorphic to the closure of
the set $ \{ w =
{z}/{(1+|z|)}\dvtx z\in{\mathbb N}^*\times{\mathbb N}^* \}$ in ${\mathbb
R}^2$. The method is
based on a ratio limit theorem for local processes and large deviation
techniques.
\end{abstract}

%

\end{frontmatter}

\section{Introduction}\label{sec0}

The concept of Martin boundary was first introduced for Brownian motion
by Martin \cite{Martin} and next extended for countable discrete time
Markov chains by Doob \cite{Doob} and Hunt \cite{Hunt}. For a Markov
chain $(Z(t))$ on a countable set $E$ with the Green function
$G(z,z')$, the Martin compactification $E_M$ is the smallest
compactification of the set $E$ for which the Martin kernels
$K(z,\cdot) = G(z,\cdot)/G(z_0,\cdot)$ extend continuously. See the
book of Woess \cite{Woess} (Chapter IV) or Rogers and Williams \cite
{Rogers} (Section III.28), for example. The Martin boundary for
homogeneous random walks in ${\mathbb Z}^d$ was obtained by Ney and
Spitzer~\cite{NeySpitzer}.

We identify the Martin boundary of a killed random walk $(Z_+(t))$ on
the positive quadrant ${\mathbb N}^*\times{\mathbb N}^*$. Such a random
walk has a substochastic transition matrix $(p(z,z') = \mu(z'-z),
z,z'\in{\mathbb N}^*\times {\mathbb N}^*)$ with some probability
measure $\mu$ on ${\mathbb Z}^2$, it is identical to a homogeneous
random walk $(S(t))$ on the two-dimensional lattice ${\mathbb Z}^2$
before it first exits from the quadrant ${\mathbb N}^*\times{\mathbb
N}^*$ and is killed at the time
\[
\tau\doteq  \inf\{n \geq0 \dvtx S(n)\notin{\mathbb N}^*\times{\mathbb
N}^*\}.
\]
The random walk $(Z_+(t))$ is therefore not homogeneous: transition
probabilities on the
boundary of the quadrant ${\mathbb N}^*\times{\mathbb N}^*$ are not
the same as in the
interior. For
nonhomogeneous Markov processes, the problem of Martin boundary
identification is usually
nontrivial and there are few examples where it was resolved.

Cartier \cite{Cartier}
described the Martin boundary of random walks on
nonhomogeneous trees and Doney \cite{Doney02} identified the Martin
boundary of
a homogeneous random walk $(Z(n))$ on ${\mathbb Z}$ killed on the
negative half-line
$\{z \dvtx z<0\}$. Alili and
Doney \cite{AliliDoney} identified the Martin boundary for space--time
random walk $S(n)
=(Z(n),n)$ for a homogeneous random walk $Z(n)$ on ${\mathbb Z}$ killed
on the
negative half-line
$\{z \dvtx z<0\}$. All these results were obtained by using a special linear
structure of the processes. The Martin boundary of Brownian motion on a
half-space was
obtained in the book of Doob \cite{Doob} by using an explicit form of the
Green function.

In Kurkova and Malyshev \cite{KurkovaMalyshev}, the full Martin
compactification is
obtained by using methods of complex analysis for nearest neighbors
random walks on a
half-plane ${\mathbb Z}\times{\mathbb N}$ and in the quadrant
${\mathbb Z}^2_+ = {\mathbb N}\times{\mathbb N}
$. In a
recent paper of
Raschel \cite{KilianRaschel}, the Martin boundary is obtained for
nearest neighbor
random walks in ${\mathbb N}\times{\mathbb N}$ with an absorption
condition on the boundary also by using methods of complex analysis.
Because of the
use of the specific algebraic setting of elliptic curves, these methods
seem to be
difficult to apply when the jump sizes are more general.

The results of Kurkova and Malyshev \cite{KurkovaMalyshev} exhibit a
formal similarity
between the limiting behavior of the Martin kernel and the optimal
large deviation
trajectories obtained by Ignatyuk, Malyshev and Scherbakov \cite
{IMS}. A natural idea is
then to study the Martin compactification by using large deviation
methods. The large
deviation approach was first proposed in the papers of
Ignatiouk-Robert \cite{Ignatiouk06,Ignatiouk08} in order to
identify the Martin
boundary for partially homogeneous random walks on a half-space
${\mathbb Z}
^{d-1}\times{\mathbb N}$. The
minimal harmonic functions were determined there by using the methods
of Choquet--Deny
theory (see Woess \cite{Woess}) and then the limiting behavior of the
Martin kernel was
obtained by using an explicit representation of the harmonic functions
combined with the
large deviation estimates of the Green function and the ratio limit
theorem of
Markov-additive processes. Unfortunately, the methods of Choquet--Deny
theory and the ratio
limit theorem are valid only for Markov-additive processes, that is, when
transition
probabilities are invariant with respect to the translations on some
directions. In the
setting of the present paper, for a random walk in the quadrant
${\mathbb N}
^*\times{\mathbb N}^*$, such an invariance property
cannot hold. Our paper is the first step toward a more ambitious
program: to identify the
Martin boundary for general partially homogeneous random walks in
${\mathbb N}^n$.

The main idea of our method is the following: to study the asymptotic
behavior of the
Martin kernels $K(z, z_n)$ for a sequence of points $z_n$ which tends
to infinity with
$\lim_n z_n/|z_n| = q$, one should consider a twisted random walk
conditioned to go to
infinity in the direction $q$. For a nonzero vector $q\in{\mathbb
R}_+^2$, such
a twisted
homogeneous random walk will visit at least one of the boundaries
$(-{\mathbb N}
)\times{\mathbb Z}$ or
${\mathbb Z}\times(-{\mathbb N})$ only a finite number of times. If
the corresponding boundary
$\{0\}\times{\mathbb N}$ (resp., ${\mathbb N}\times\{0\}$) is removed,
the resulting
process is then
identical to the homogeneous random walk $(S(t))$ before the first time
when it hits the
set ${\mathbb Z}\times(-{\mathbb N})$ [resp., $(-{\mathbb N})\times
{\mathbb Z}$]. The limiting
behavior of
the Martin kernel
of this process corresponding to the direction $q$ is already known in
such a setting. The
limiting behavior of the Martin kernel of the original process
$(Z_+(t))$ should be
essentially the same but with a correction given by a potential
function. When both
coordinates of $q$ are positive, this idea is transformed into a
rigorous proof with the
aid of large deviation estimates and a generalization of a ratio limit
theorem of the
paper \cite{Ignatiouk06}. When one of the coordinates of $q$ is zero,
that is, when the
process is conditioned to go to infinity along one of the boundaries,
our proof is much more
complicated. In this case, we combine large deviation techniques and
the ratio limit theorem with
delicate estimates obtained from the Harnack inequalities.

We assume that the probability measure $\mu$ on ${\mathbb Z}^2$ satisfies
the following conditions:
\begin{itemize}[(H3)]
\item[(H1)]\hypertarget{hypoH1} \mbox{\textit{The homogeneous random
walk}} $S(t)=(S_1(t),S_2(t))$ \textit{on} ${\mathbb Z}^2$ \textit{having transition
probabilities}
$p_S(z,z')=\mu(z'-z)$ \textit{is irreducible and}
\[
m \doteq  \sum_{z\in{\mathbb Z}^d} z \mu(z) \not= 0.
\]
\item[(H2)]\hypertarget{hypoH2}\mbox{ }\textit{The killed random walk} $(Z_+(t))$ \textit{is
irreducible on} ${\mathbb N}^*\times{\mathbb N}^*$.
\item[(H3)]\hypertarget{hypoH3}\mbox{ }\textit{The jump generating
function}
%
%
\begin{equation}\label{e1-1}
\varphi(a) \doteq  \sum_{z\in{\mathbb Z}^2} \mu(z)
\exp(a\cdot z)
\end{equation}
\textit{is finite everywhere on} ${\mathbb R}^2$.
\item[(H4)]\hypertarget{hypoH4} $(S_1(t))$ \textit{and} $(S_2(t))$ \textit{are
aperiodic random walks on} ${\mathbb Z}$.
\end{itemize}
Under the above assumptions, the set
\[
D\doteq  \{a\in{\mathbb R}^2 \dvtx\varphi(a)\leq1\}
\]
is compact and strictly convex, the gradient $\nabla\varphi(a)$ exists
everywhere on
${\mathbb R}^2$ and does not vanish on the
boundary $\partial D = \{a\in{\mathbb R}^2 \dvtx\varphi(a)= 1\}$, the mapping
%
%
\begin{equation}\label{e1-2}
a\to q(a) \doteq  \nabla\varphi(a)/|\nabla\varphi(a)|
\end{equation}
determines a homeomorphism from $\partial D$ to the unit
two-dimensional sphere ${\mathcal
S}^2=\{q\in{\mathbb R}^2\dvtx |q| = 1\}$ (see \cite{Hennequin}). We
denote by
$q\to
a(q)$ the inverse mapping of
(\ref{e1-2}) and we let $a(q) = a(q/|q|)$ for a nonzero $q\in{\mathbb R}^2$.
According to this
notation, $a(q)$ is the only point in $\partial D$ where the
vector $q$ is normal to the convex set $D$. Throughout this paper, we
denote by ${\mathbb N}$ the
set of all nonnegative integers and we let ${\mathbb N}^*={\mathbb
N}\setminus\{0\}$. The
set of all
nonnegative real numbers is denoted by ${\mathbb R}_+ =[0, +\infty[$
and ${\mathbb R}^*_+
= \ ]0,+\infty[$
denotes the set of all strictly positive real numbers.
It is convenient moreover to introduce the following notation:
${\mathbb N}$
denotes the set
of all nonnegative integers and ${\mathbb N}^* ={\mathbb N}\setminus\{
0\}$,
\[
\tau\doteq  \inf\{n \geq0\dvtx S(n)\notin{\mathbb N}^*\times{\mathbb
N}^*\}
\]
is the first time when the random walk $(S(t))$ exits from the quadrant
${\mathbb N}^*\times{\mathbb N}^*$,
\[
{\mathcal
S}^2_+\doteq \{q\in{\mathbb R}_+^2\dvtx |q| = 1\} \quad\mbox{and}\quad
\Gamma_+ \doteq  \{a\in\partial D \dvtx q(a)\in{\mathcal S}^2_+\}.
\]
For $a\in\Gamma_+$ and $z=(x_1,x_2)\in{\mathbb N}^*\times{\mathbb
N}^*$, we set
%
%
\begin{equation}\label{e1-3}\qquad
h_a(z) \doteq  \cases{
x_1\exp(a\cdot z) - {\mathbb E}_z\bigl( S_1(\tau)\exp\bigl(a\cdot S(\tau)\bigr),
\tau<
\infty
\bigr), \cr \hspace*{174.43pt}\hspace*{-0.92pt}\qquad\mbox{if
$q(a)=(0,1)$}, \cr
x_2\exp(a\cdot z) - {\mathbb E}_z\bigl( S_2(\tau)\exp\bigl(a\cdot S(\tau)\bigr),
\tau<
\infty
\bigr), \cr
\qquad\hspace*{174.43pt}\hspace*{-0.92pt}\mbox{if
$q(a)=(1,0)$},\cr
\exp(a\cdot z) - {\mathbb E}_z\bigl( \exp\bigl(a\cdot S(\tau)\bigr), \tau< \infty
\bigr), \qquad\mbox{otherwise.}}
\end{equation}
$G_+(z,z')$ denotes the Green function of the process $(Z(t))$:
\[
G_+(z,z') = \sum_{n=0}^\infty{\mathbb P}_z\bigl(Z_+(n)=z'\bigr).
\]

The main result of our paper is the following theorem.
\begin{theorem}\label{th1} Under the hypotheses \textup{\hyperlink{hypoH1}{(H1)}--\hyperlink{hypoH4}{(H4)}}, for any
$q\in{\mathcal S}^2_+$ and any sequence of points
$z_n\in{\mathbb N}^*\times{\mathbb N}^*$ with $\lim_n|z_n| = \infty
$ and $\lim_n
z_n/|z_n| = q$,
%
%
\begin{equation}\label{e1-4}
\lim_{n\to\infty} {G_+(z,z_n)}/{G_+(z_0,z_n)} = h_{a(q)}(z)/h_{a(q)}(z_0)
\end{equation}
for all $z\in{\mathbb N}^*\times{\mathbb N}^*$.
\end{theorem}

Remark that the conditions \hyperlink{hypoH1}{(H1)} and
\hyperlink{hypoH2}{(H2)} are essential for our approach, our method
does not work when at least one of them is not satisfied. The
hypotheses \hyperlink{hypoH3}{(H3)} and \hyperlink{hypoH4}{(H4)} are
required by the paper \cite{Ignatiouk06}, we use its results to get
(\ref{e1-4}) for $q\in\{(1,0), (0,1)\}$. When the coordinates $q_1$ and
$q_2$ of the vector $q = \lim_n z_n/|z_n|$ are nonzero, the assumption
\hyperlink{hypoH4}{(H4)} is not needed and the hypotheses
\hyperlink{hypoH3}{(H3)} can be replaced by a less restrictive
condition of Ney and Spitzer \cite {NeySpitzer} where the jump
generating function (\ref{e1-1}) is assumed to be finite only in a
neighborhood of the set~$D$.

Recall that a sequence $z_n$ is said to converge to a point on the
Martin boundary $\partial_M({\mathbb N}^*\times{\mathbb N}^*)$ of
${\mathbb N}^*\times{\mathbb N}^*$ determined by the Markov process
$(Z_+(t))$ if and only if the sequence of functions $
z\to{G_+(z,z_n)}/{G_+(z_0,z_n)}$ converges point-wise on ${\mathbb
N}^*\times{\mathbb N}^*$. According to this definition, Theorem
\ref{th1} implies the following statement.
\begin{cor}\label{cor1.1}
Under the hypotheses
\textup{\hyperlink{hypoH1}{(H1)}--\hyperlink{hypoH4}{(H4)}}, the following
assertions hold:

(1) A sequence of points
$z_n\in{\mathbb N}^*\times{\mathbb N}^*$ with $\lim_n |z_n| =
+\infty$ converge to a point of the Martin boundary for the Markov
process $Z_+(t)$ if and
only if $z_n/|z_n| \to q$ for
some point $q\in{\mathcal S}^2_+$.

(2) The full Martin compactification of the quadrant
${\mathbb N}^*\times{\mathbb N}^*$ is homeomorphic to the closure of
the set $ \{ w =
{z}/{(1+|z|)} \dvtx z\in{\mathbb N}^*\times{\mathbb N}^* \}$ in ${\mathbb R}^2$.
\end{cor}

Our paper is organized as follows. In Section \ref{sec3}, the main idea
of the proof
of our result is sketched. Section \ref{sec2} is devoted to the
preliminary results. In Section \ref{Sec4}, we prove that the functions
$h_a$ with
$a\in\Gamma_+$ defined by (\ref{e1-3}) are finite, harmonic for the
Markov process
$(Z_+(t))$ and strictly positive. Section \ref{sec5} is devoted to the
large deviation
results. It is shown that the family of scaled processes
$Z_+^\varepsilon(t) =
\varepsilon Z_+([t/\varepsilon])$
satisfies sample path large deviation principle. The logarithmic
estimates of the Green
function are obtained from the corresponding large deviation bounds. In
Section \ref{sec6}, the large deviation estimates are used to decompose
the Green function
$G_+(z,z_n)$ into a main part corresponding to an optimal large
deviation way to go from
$z$ to $z_n$ and the negligible part. In Section \ref{sec7}, we
generalize the ratio
limit theorem of Ignatiouk-Robert \cite{Ignatiouk06}. The
decomposition into a main and a
negligible parts of the Green function $G_+(z,z_n)$ and the ratio limit
theorem are
next combined in Section \ref{sec8} in order to complete the proof of
Theorem \ref{th1}.

\section{Local processes and renewal equations: A sketch of
proofs}\label{sec3}
The main steps of our method can be summarized as follows:

(1) For a sequence $(z_n)\in{\mathbb N}^*\times{\mathbb N}^*$ with
$\lim_n z_n/|z_n| = q$ and $\lim_n|z_n| = +\infty$, the Green function
$G_+(z,z_n)$ of the Markov process $(Z_+(t))$ is represented in terms
of a local random walk which is Markov-additive and has the same
transition probabilities as the original random walk $(Z_+(t))$ in a
neighborhood of the point $q|z_n|$.

(2) Next, large deviation estimates are used to decompose $G_+(z,z_n)$
into a main part corresponding to an optimal large deviation way to go
from $z$ to $z_n$ and the negligible part. Such a decomposition allows
us to get the limit of the Martin kernel
\[
\lim_n G_+(z,z_n)/G_+(z_0,z_n)
\]
from the limiting behavior and the uniform bounds of the Martin kernel
of the corresponding local process.

When the coordinates of the vector $q=(q_1,q_2)$ are nonzero, the local
Markov-additive process is simply a homogeneous random walk $(S(t))$ on
${\mathbb Z}^2$ having
transition probabilities $p(z,z') = \mu(z'-z)$. This is the simplest
case in our proof.
The following renewal equation represents the Green
function $G_+(z,z')$ of the Markov process $(Z_+(t))$ in terms of the
Green function
$G(z,z')$ of the random walk $(S(t))$:
%
%
\begin{equation}\label{e2-1}
G_+(z,z') = G(z,z') - {\mathbb E}_z \bigl(G(S(\tau),z'),
\tau<\infty\bigr).
\end{equation}
Ney and Spitzer \cite{NeySpitzer} proved
that for any $q\in{\mathcal S}^2_+$ and any sequence of points
$z_n\in{\mathbb N}^*\times{\mathbb N}^*$ with ${\lim_n}|z_n| = \infty
$ and $\lim_n
z_n/|z_n| = q$,
%
%
\begin{equation}\label{e2-2}
\lim_{n\to\infty} G(z,z_n)/G(0,z_n) = \exp\bigl(a(q)\cdot z\bigr)
\end{equation}
for all $z\in{\mathbb Z}^2$ (see also Section 7 in \cite
{Ignatiouk06} for an
alternative simple proof of this result). Using the renewal equation
(\ref{e2-1}), one can therefore get the
equality
%
%
\begin{eqnarray}\label{e2-3}
\lim_{n\to\infty} \frac{G_+(z,z_n)}{G(0,z_n)} &=& \exp\bigl(a(q)\cdot z\bigr) -
{\mathbb E}_z \bigl(\exp\bigl(a(q)\cdot S(\tau)\bigr), \tau<\infty\bigr)
\nonumber\\[-8pt]\\[-8pt]
&\doteq & h_{a(q)}(z),\nonumber
\end{eqnarray}
if one can prove the exchange of limits
%
%
\begin{equation}\label{e2-4}\quad
\lim_{n\to\infty} {\mathbb E}_z \biggl(\frac{G(S(\tau) ,z_n)}{G(0,z_n)},
\tau
<\infty\biggr)
= {\mathbb E}_z \biggl(\lim_{n\to\infty}\frac{G(S(\tau)
,z_n)}{G(0,z_n)}, \tau
<\infty\biggr).
\end{equation}
Relation (\ref{e1-4}) will follow finally from the relation (\ref
{e2-3}) because the function $h_a$ is strictly positive on ${\mathbb
N}^*\times{\mathbb N}^*$ (see Proposition \ref{Pr4-1} below). Equality
(\ref{e2-4}) is therefore a key relation for our problem.

While the above idea seems quite simple, the proof of (\ref{e2-4}) is
nontrivial because the convergence (\ref{e2-2}) is not uniform and the
classical convergence theorems are here difficult to use. With our
approach, for a sequence of points $z_n\in{\mathbb N}^*\times{\mathbb
N}^*$ with ${\lim_n}|z_n|= \infty$ and $\lim_n z_n/|z_n| = q$, we first
decompose the right-hand side of (\ref{e2-1}) into a main part
\[
\Xi^q_\delta(z,z_n) \doteq  G(z,z_n) - {\mathbb E}_z
\bigl(G(S(\tau),z_n), \tau <\infty , |S(\tau)| < \delta|z_n| \bigr)
\]
and the corresponding negligible part by using the large deviation
estimates of the Green function $G(z,z')$ and $G_+(z,z')$. Next, we get
the estimates
%
%
\begin{equation}\label{e2-5}
\sup_n \vienas_{\{|z| <\delta|z_n|\}} G(z,z_n)/G(z_0,z_n) \leq C(z)
\end{equation}
such that $ {\mathbb E}_z(C(S(\tau)), \tau<\infty) < \infty$ and
finally, using the point-wise convergence (\ref{e2-2}) and dominated
convergence theorem we obtain (\ref{e1-4}). The estimates (\ref{e2-5})
are obtained in Section \ref{sec7} with a suitable exponential function
$C(z)$ by using the ratio limit theorem applied to the random walk
$(S(t))$.

The case when one of the coordinates of the vector $q$ is equal to
zero, that is, when the sequence $(z_n)$ tends to infinity along one of
the boundaries of the domain, is much more delicate to handle. First of
all, we cannot use here the renewal equation (\ref{e2-1}) because the
function $ \exp(a(q)\cdot z) - {\mathbb E}_z( \exp(a(q)\cdot S(\tau)),
\tau <\infty) $ is in this case identical to zero. If $q=(1,0)$, one
should consider a Markov-additive process having the same statistical
behavior as the process $(Z_+(t))$ near the boundary ${\mathbb N}
\times\{0\}$ and far from the boundary $\{0\}\times{\mathbb N}$. This
is a random walk $(Z^1_+(t))$ on ${\mathbb Z}\times{\mathbb N}^*$
having a substochastic transition matrix $(p_1(z,z') = \mu(z'-z),
z,z'\in{\mathbb Z}\times{\mathbb N}^*)$. It is identical to the random
walk $(S(t))$ before the time $ \tau_2 \doteq  \inf\{t \geq0 \dvtx
S_2(t)\leq0\}$ and killed at the time $\tau_2$. Our Markov process
$(Z_+(t))$ is therefore identical to $(Z^1_+(t))$ before the time $
\tau_1 \doteq  \inf\{t \geq0 \dvtx S_1(t)\leq0\}$. Since clearly $\tau=
\min\{\tau_1, \tau_2\}$, the Green function $G_+(z,z')$ of the Markov
process $(Z_+(t))$ is related to the Green function $G^1_+(z,z')$ of
the process $(Z^1_+(t))$ as follows:
%
%
\begin{equation}\label{e2-6}
G_+(z,z') = G^1_+(z,z') - {\mathbb E}_z \bigl(G^1_+(S(\tau),z'), \tau=\tau_1
<\tau_2 \bigr).
\end{equation}
Theorem 1 of \cite{Ignatiouk06} proves that for any sequence of points
$z_n\in{\mathbb N}^*\times{\mathbb N}^*$ with ${\lim_n}|z_n| = \infty $
and $\lim_n z_n/|z_n| = q=(1,0)$,
%
%
\begin{equation}\label{e2-7}\hspace*{28pt}
\lim_{n\to\infty} G^1_+(z,z_n)/G^1_+(z_0,z_n) =
h_{a(q),+}^1(z)/h_{a(q),+}^1(z_0)\qquad
\forall z\in{\mathbb Z}\times{\mathbb N}^*,
\end{equation}
with a strictly positive function $h_{a(q),+}^1$ on ${\mathbb Z}\times
{\mathbb N}^*$ defined by
\[
h_{a(q),+}^1(z) = x_2\exp\bigl(a(q)\cdot z\bigr) -
{\mathbb E}_z \bigl(S_2(\tau_2)\exp\bigl(a(q)\cdot S(\tau_2)\bigr), \tau_2 <\infty\bigr).
\]
Similarly to the previous case, we decompose the right-hand side of
the renewal equation (\ref{e2-6}) into a main part
\[
G^1_+(z,z_n) -
{\mathbb E}_z \bigl(G^1_+(S(\tau), z_n), \tau= \tau_1 < \tau_2, |S(\tau
)|<\delta
|z_n| \bigr)
\]
and the corresponding negligible part
by using the large deviation estimates of the Green functions
$G_+(z,z')$ and
$G^1_+(z,z')$ and we show
there are $\delta> 0$ and a function $C_+^1(z)$ with
\[
{\mathbb E}_z\bigl(C_+^1(S(\tau)), \tau= \tau_1 <\tau_2\bigr) < \infty
\]
such that
%
%
\begin{equation}\label{e2-8}
\sup_n \vienas_{\{|z| <\delta|z_n|\}} G_+^1(z,z_n)/G_+^1(z_0,z_n)
\leq C_+^1(z).
\end{equation}
The proof of these estimates is the most delicate part of our work.

\section{Preliminary results}\label{sec2}
For a given $a\in D \doteq  \{a\in{\mathbb R}^2 \dvtx \varphi(a)\leq1\}$,
let us
consider a new twisted homogeneous random walk $(S^a(t))$ on
${\mathbb Z}^2$ having transition probabilities
%
%
\begin{equation}\label{e3-1}
p_a(z,z') = \mu(z'-z)\exp\bigl(a\cdot(z'-z)\bigr).
\end{equation}
According to the
definition of the set $D$, the transition matrix of such a random walk is
substochastic. Recall that
\[
\tau= \tau_1\wedge\tau_2,
\]
where
$\tau_1 \doteq  \inf\{n \geq0 \dvtx S(n)\notin{\mathbb N}^*\times
{\mathbb Z}\}$ and $
\tau_2 \doteq  \inf\{n \geq0 \dvtx S(n)\notin{\mathbb Z}\times
{\mathbb N}^*\}$.
\begin{prop}\label{Pr3-1} For every $a\in D$, the quantity $
{\mathbb E}_z( \exp(a\cdot(S(\tau)-z)), \tau< \infty)$
is equal to the probability that the twisted random walk $(S^a(t))$
starting at the point
$z$ ever exits from the positive quadrant ${\mathbb N}^*\times{\mathbb N}^*$.
\end{prop}
\begin{pf} Indeed, let $
\tau^a $ denote the first time when the twisted random walk $(S^a(t))$
exits from the quadrant
${\mathbb N}^*\times{\mathbb N}^*$. Then for any $t\in{\mathbb N}$,
\[
{\mathbb P}_z\bigl(S^a(t) = z', \tau=t\bigr) = \exp\bigl(a\cdot(z'-z)\bigr) {\mathbb
P}_z\bigl(S(t)=z', \tau=t\bigr)\qquad
\forall z,z'\in{\mathbb Z}^2
\]
and consequently, ${\mathbb P}_z(\tau^a < \infty) = {\mathbb E}_z(
\exp(a\cdot
(S(\tau
)-z)), \tau< \infty)$.
\end{pf}

The
set $\Gamma_+= \{a\in\partial D \dvtx q(a)\in{\mathcal S}_+^2\}$ endowed
with a topology
induced by the usual topology of ${\mathbb R}^2$ is homeomorphic to a segment
with the end points
in $a(1,0)$ and $a(0,1)$. The points $a(1,0)$
and $a(0,1)$ are said to be critical.
\begin{prop}\label{Pr3-2} Every noncritical point
of $\Gamma_+$ has a neighborhood where the functions $
a \to{\mathbb E}_z(\exp(a\cdot S(\tau)), \tau< \infty)$
are finite for all $z\in{\mathbb N}^*\times{\mathbb N}^*$.
\end{prop}
\begin{pf} By
Proposition \ref{Pr3-1}, the
function $a \to{\mathbb E}_z(\exp(a\cdot S(\tau)), \tau< \infty)$ is
finite on
$D\doteq \{a\in{\mathbb R}^2 \dvtx \varphi(a)\leq1\}$.
Furthermore, let us consider the critical points $a(1,0)=(a'_1,a'_2)$ and
$a(0,1)=(a_1'',a_2'')$. Recall that under the hypotheses \hyperlink{hypoH1}{(H1)} and \hyperlink{hypoH3}{(H3)}
the set $D$ is
compact and strictly
convex, and according to the definition of the mapping $q\to a(q)$,
\[
\nabla\varphi(a'_1,a'_2) = |\nabla\varphi(a'_1,a'_2)| (1,0)
\quad\mbox{and}\quad
\nabla\varphi(a''_1,a''_2) = |\nabla\varphi(a''_1,a''_2)| (0,1).
\]
Every noncritical point of $\Gamma_+$ has therefore a neighborhood
where for any point
$a=(a_1,a_2)\notin D$
there exist two points $\hat{a} = (\hat{a}_1,\hat{a}_2)$ and $\tilde
{a}=(\tilde{a}_1,\tilde{a}_2)$ on
the boundary of the set $D$ with $\hat{a}_1=a_1$, $\hat{a}_2 < a_2$ and
$\tilde{a}_1
< a_1$, $\tilde{a}_2 = a_2$
%
%
\begin{figure}

\includegraphics{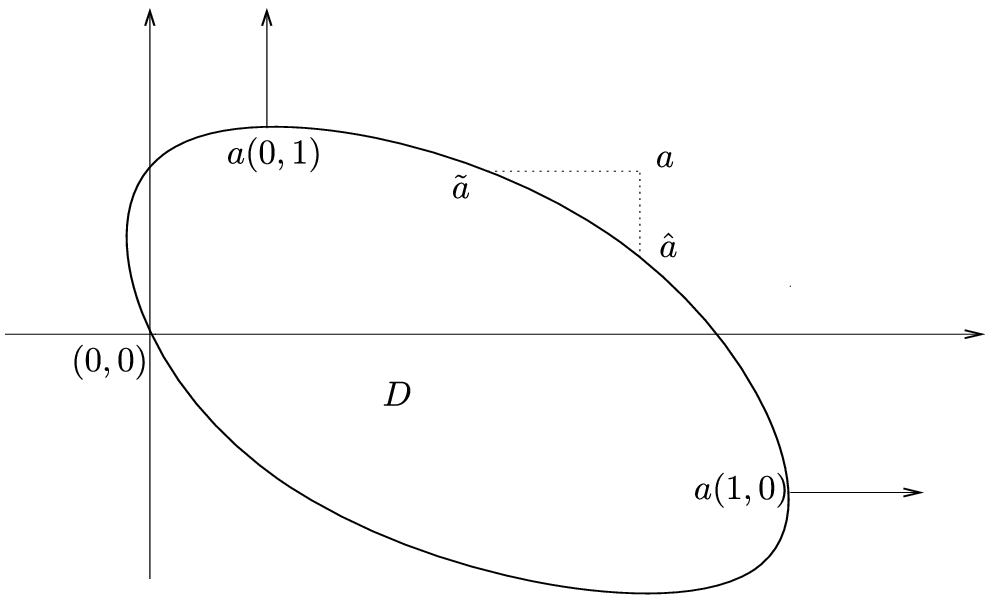}

\caption{ }
\label{figure1}
\end{figure}
(see Figure~\ref{figure1}).
Since $S_1(\tau)\leq0$ on the event $\{\tau=\tau_1 < +\infty\}$, and
$S_2(\tau)\leq
0$ on the event $\{\tau=\tau_2 < +\infty\}$, from this it follows that
\begin{eqnarray*}
&&{\mathbb E}_z\bigl(\exp\bigl(a\cdot S(\tau)\bigr), \tau< +\infty\bigr) \\
&&\qquad\leq
{\mathbb E}_z\bigl(\exp
\bigl(a\cdot
S(\tau)\bigr), \tau=\tau_1 <
+\infty\bigr) + {\mathbb E}_z\bigl(\exp\bigl(a\cdot S(\tau)\bigr), \tau=\tau_2 <
+\infty\bigr)
\\
&&\qquad\leq{\mathbb E}_z\bigl(\exp\bigl(\tilde{a}\cdot S(\tau)\bigr), \tau= \tau_1 <
+\infty\bigr) + {\mathbb E}_z\bigl(\exp\bigl(\hat{a}\cdot S(\tau)\bigr), \tau= \tau_2 <
+\infty\bigr)
\\
&&\qquad\leq{\mathbb E}_z\bigl(\exp\bigl(\tilde{a}\cdot S(\tau)\bigr), \tau<
+\infty\bigr) + {\mathbb E}_z\bigl(\exp\bigl(\hat{a}\cdot S(\tau)\bigr), \tau< +\infty\bigr)
< +\infty.\qquad
\end{eqnarray*}
\upqed\end{pf}
\begin{prop}\label{Pr3-3} The critical point
$a(1,0)=(a_1',a_2')$ has a neighborhood where the functions $a\to
{\mathbb E}
_z(\exp(a\cdot S(\tau)), \tau
= \tau_1 < \tau_2)$ are finite for all \mbox{$z\in{\mathbb N}^*\times
{\mathbb N}^*$}. Moreover,
for any $\delta>
0$ small enough there is a point $\hat{a} =(\hat{a}_1,\hat{a}_2)\in
\partial D$ with
$\hat{a}_1 < a_1'$ and $\hat{a}_2 = a'_2 + \delta$ such that
%
%
\begin{equation}\label{e3-2}
{\mathbb E}_z\bigl(\exp\bigl(a(1,0)\cdot S(\tau) + \delta S_2(\tau)\bigr), \tau=
\tau_1 < \tau _2\bigr) \leq \exp(\hat{a}\cdot z)
\end{equation}
for all $z\in{\mathbb N}^*\times{\mathbb N}^*$.
\end{prop}
\begin{pf} The proof of this proposition uses essentially the same
arguments as the
proof of Proposition \ref{Pr3-2}. For $a\in D$,
\[
{\mathbb E}_z\bigl(\exp\bigl(a\cdot S(\tau)\bigr), \tau= \tau_1 < \tau_2\bigr) \leq
\exp(a\cdot z)\qquad \forall z\in{\mathbb N}^*\times{\mathbb N}^*,
\]
because the quantity ${\mathbb E}_z(\exp(a\cdot(S(\tau)-z)), \tau=
\tau_1 < \tau_2)$ is equal to the probability that the twisted
substochastic homogeneous random walk $(S^a(t))$ starting at $z$ hits
the set $(-{\mathbb N})\times{\mathbb Z}$ before hitting the set
${\mathbb Z}\times(-{\mathbb N})$. This proves that the functions
$a\to{\mathbb E}_z(\exp (a\cdot S(\tau)), \tau = \tau_1 < \tau_2)$ are
finite on $D$ for all $z\in{\mathbb N}^*\times{\mathbb N}^*$. Moreover,
let us consider the points $a(0,1)=(a_1'',a_2'')$ and
$a(0,-1)=(a_1''',a_2''')$ on the boundary $\partial D$ of $D$. Then the
set $ \Omega\doteq  \{a=(a_1,a_2) \in{\mathbb R}^2 \dvtx a_1 >
\max\{a_1'', a_1'''\}, a_2''' < a_2 < a_2''\}$ is an open neighborhood
of the point $a(1,0)$ and for any $a=(a_1,a_2)\in \Omega\setminus D$
there is a point $\hat{a} =(\hat{a}_1, \hat{a}_2)$ on the boundary of
the set $D$ with $\hat{a}_2=a_2$ and $\hat{a}_1 < a_1$ (see Figure
\ref{figure2}).
%
%
\begin{figure}

\includegraphics{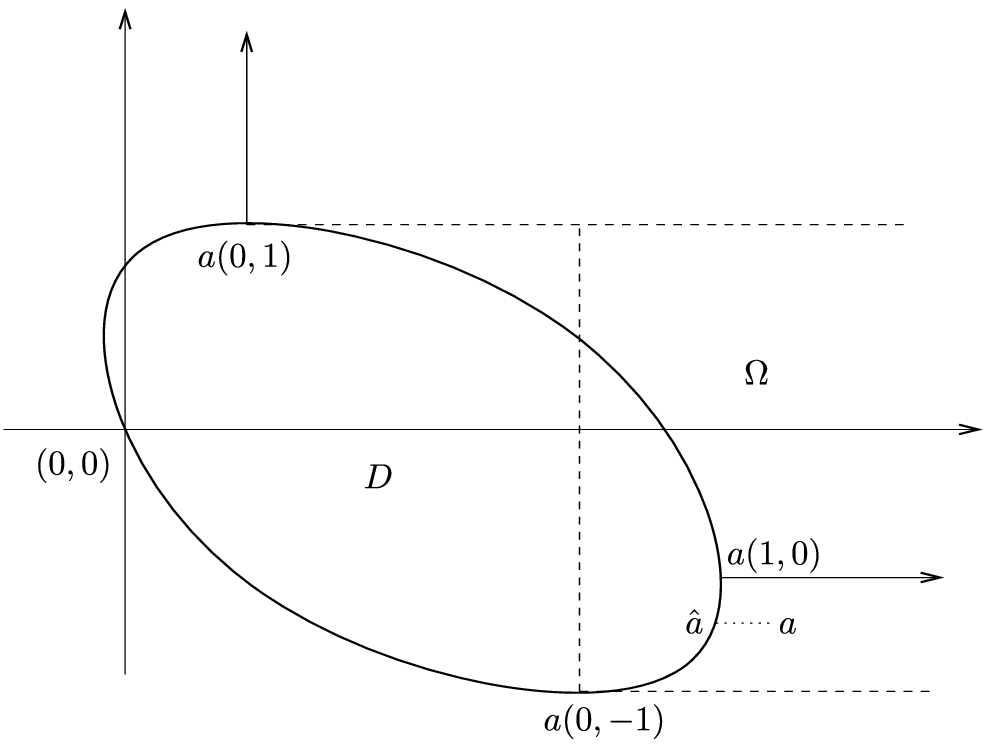}

\caption{ }
\label{figure2}
\end{figure}
Since on the event $\{ \tau=\tau_1 <
\tau_2\}$, $S_1(\tau) \leq0$ we conclude that for any $z\in{\mathbb N}
^*\times{\mathbb N}^*$,
\[
{\mathbb E}_z\bigl(\exp\bigl(a\cdot S(\tau)\bigr), \tau=\tau_1 < \tau_2\bigr) \leq
{\mathbb E}_z\bigl(\exp
\bigl(\hat
{a}\cdot
S(\tau)\bigr), \tau=\tau_1
< \tau_2\bigr) \leq\exp(\hat{a}\cdot z).
\]
The functions $a\to{\mathbb E}_z(\exp(a\cdot S(\tau)), \tau
= \tau_1 < \tau_2)$ are therefore finite on $\Omega$ for all
$z\in{\mathbb N}^*\times{\mathbb N}^*$. Finally, for $\delta>0$
small enough, $a=a(1,0)
+ (0,\delta)\in
\Omega$ and hence, the last inequality proves also (\ref{e3-2}).
\end{pf}

A straightforward consequence of Proposition \ref{Pr3-3} is the
following statement.
\begin{cor}\label{Cor3-1} For $a=a(1,0)$, the function
%
%
\begin{equation}\label{e3-3}
z\to{\mathbb E}_z\bigl(|S_2(\tau)|\exp\bigl(a\cdot S(\tau)\bigr), \tau
= \tau_1 < \tau_2\bigr)
\end{equation}
is finite on ${\mathbb N}^*\times{\mathbb N}^*$.
\end{cor}
\begin{pf} Indeed, on the event $\tau
= \tau_1 < \tau_2$, for any $\delta> 0$, one has
\[
0 < S_2(\tau) \leq\frac{1}{\delta} \exp(\delta S_2(\tau))
\]
and consequently, for $a=a(1,0)$,
\begin{eqnarray*}
&&{\mathbb E}_z\bigl(|S_2(\tau)|\exp\bigl(a\cdot S(\tau)\bigr), \tau=\tau_1< \tau
_2\bigr)\\
&&\qquad=
{\mathbb E}_z\bigl(S_2(\tau)\exp\bigl(a\cdot S(\tau)\bigr), \tau=\tau_1< \tau_2\bigr)
\\
&&\qquad\leq\frac{1}{\delta} {\mathbb E}_z\bigl(\exp\bigl(a\cdot S(\tau) + \delta
S_2(\tau)\bigr),
\tau=\tau_1< \tau_2\bigr).
\end{eqnarray*}
Since by Proposition \ref{Pr3-3}, the right-hand side of the last relation
is finite for all $z\in{\mathbb N}^*\times{\mathbb N}^*$ and $\delta
> 0$ small enough,
we conclude that
the function (\ref{e3-3}) is finite on ${\mathbb N}^*\times{\mathbb N}^*$.
\end{pf}

To show that the functions (\ref{e1-3}) are well defined, we will need,
moreover, the following statement.
\begin{lemma}\label{Lem3-1} For a random walk $(\xi(t))$ on ${\mathbb
Z}$ having
zero mean and transition
probabilities $P(x,x')=P(0,x'-x)$ such that for some $\delta> 0$,
\[
\sum_x e^{-\delta x} P(0,x) < \infty\quad\mbox{and}\quad\sum_x |x| P(0,x) <
\infty,
\]
the function $f(x) = {\mathbb E}_x(|\xi(T_0)|)$
with $T_0 = \inf\{t\geq0 \dvtx \xi(t)\leq0\}$ is finite everywhere on
${\mathbb N}^*$.
\end{lemma}

This elementary lemma has been proved in the proof of Lemma 5.3 in
Ignatiouk~\cite{Ignatiouk06}. A more general related result can also be
found in Chow \cite{Chow}. Corollary~\ref{Cor3-1} combined with Lemma
\ref{Lem3-1} implies the following proposition.
\begin{prop}\label{Pr3-4}
The function $z \to{\mathbb E}_z(|S_2(\tau)|\exp (a(1,0)\cdot S(\tau)),
\tau< \infty)$ is finite on ${\mathbb N}^*\times{\mathbb N}^*$.
\end{prop}
\begin{pf} To prove this proposition, let us first notice that
\begin{eqnarray*}
&&{\mathbb E}_z\bigl(|S_2(\tau)|\exp\bigl(a\cdot S(\tau)\bigr), \tau< \infty\bigr) \\
&&\qquad=
{\mathbb E}_z\bigl(|S_2(\tau)|\exp\bigl(a\cdot S(\tau)\bigr), \tau=\tau_1< \tau
_2\bigr) \\
&&\qquad\quad{} + {\mathbb E}_z\bigl(|S_2(\tau)|\exp\bigl(a\cdot
S(\tau)\bigr), \tau=\tau_2< \infty\bigr),
\end{eqnarray*}
where for $a=a(1,0)$, by Corollary \ref{Cor3-1},
\[
{\mathbb E}_z\bigl(|S_2(\tau)|\exp\bigl(a(1,0)\cdot S(\tau)\bigr), \tau=\tau_1<
\tau_2\bigr) <
\infty\qquad
\forall z\in{\mathbb N}^*\times{\mathbb N}^*.
\]
To prove that the function $z \to{\mathbb E}_z(|S_2(\tau)|\exp
(a(1,0)\cdot
S(\tau)), \tau< \infty)$ is finite on ${\mathbb N}^*\times{\mathbb
N}^*$ it is
therefore sufficient to show that
%
%
\begin{equation}\label{e3-4}\quad
{\mathbb E}_z\bigl(|S_2(\tau)|\exp\bigl(a(1,0)\cdot S(\tau)\bigr), \tau=\tau_2<
\infty\bigr) <
\infty\qquad
\forall z\in{\mathbb N}^*\times{\mathbb N}^*.
\end{equation}
Next, we consider a twisted random walk $(S^a(t))$ on ${\mathbb Z}^2$
with transition
probabilities $p_a(z,z') = \mu(z'-z)\exp(a\cdot(z'-z))$ for $a=a(1,0)$.
The second coordinate $(S_2^a(t))$ of $(S^a(t))$ is a random walk on
${\mathbb Z}$ having a mean
\[
{\mathbb E}_0(S_2^a(1)) = \frac{\partial}{\partial a_2}\varphi(a_1,a_2)
\bigg|_{(a_1,a_2)=a(1,0)}
= 0
\]
and satisfying the conditions of Lemma \ref{Lem3-1}. This lemma applied
with $\xi(t)=S_2^a(t)$ and $T_0= \tau_2^a \doteq  \inf\{n \geq0 \dvtx
S^a_2(n) \leq0\}$ proves that the function
${\mathbb E}_x(|S_2^a(\tau^a_2)|)$ is finite on ${\mathbb N}^*$. Since
for any
$z=(x_1,x_2)\in{\mathbb N}^*\times{\mathbb N}^*$,
\[
{\mathbb E}_z\bigl(|S_2(\tau_2)|\exp\bigl(a(1,0)\cdot
S(\tau_2)\bigr), \tau_2 < \infty\bigr) = {\mathbb E}_{x_2}(|S_2^a(\tau_2^a)|)
\]
we conclude that (\ref{e3-4}) holds. Proposition \ref{Pr3-4} is
therefore proved.
\end{pf}

\section{Harmonic functions}\label{Sec4}

The main result of this section is the following proposition.
\begin{prop}\label{Pr4-1} For every $a\in\Gamma_+$, the functions $h_a$
defined by
(\ref{e1-3}) is finite, strictly positive on ${\mathbb N}^*\times
{\mathbb N}^*$ and
harmonic for the Markov process $(Z_+(t))$.
\end{prop}

Before proving this proposition, we consider the following lemmas.
\begin{lemma}\label{Lem4-1} For $a\in\Gamma_+$, the function
$
z \to1 - {\mathbb E}_z ( \exp(a\cdot(S(\tau)-z)), \tau< \infty)
$
is strictly positive on ${\mathbb N}^*\times{\mathbb N}^*$ when
$q(a)\notin\{(1,0),
(0,1)\}$ and is
identically zero when $q(a)\in\{(1,0), (0,1)\}$.
\end{lemma}
\begin{pf} Indeed, for any $a\in\Gamma_+$,
the twisted random walk $S^a(t)=(S_1^a(t),\break S_2^2(t))$ has a stochastic
transition matrix $(p_a(z,z') =
\exp(a\cdot(z'-z))\mu(z'-z),\break z,z'\in{\mathbb Z}^2)$, a nonzero mean
\[
m(a) = \sum_{z\in{\mathbb Z}^2} z \exp(a\cdot z) \mu(z) = \nabla
\varphi
(a) =
|\nabla\varphi(a)| q(a)
\]
and a finite variance. If $q(a) = (0,1)$, the first coordinate
$S^a_1(t)$ of $S^a(t)$ is
therefore a recurrent random walk on ${\mathbb Z}$, the first time when
$S^a_1(t)$
becomes negative or zero is almost surely finite for any starting point
$S^a(0)=z\in{\mathbb N}^*\times{\mathbb N}^*$ and consequently, the
twisted random walk
$(S^a(t))$ almost
surely exits from the quadrant ${\mathbb N}^*\times{\mathbb N}^*$. By
Proposition \ref
{Pr3-1}, from this it
follows that
\[
1 - {\mathbb E}_z\bigl ( \exp\bigl(a\cdot\bigl(S(\tau)-z\bigr)\bigr), \tau< \infty\bigr) =
{\mathbb P}_z(\tau^a = \infty) = 0\qquad \forall z\in{\mathbb N}^*\times
{\mathbb N}^*.
\]
The same arguments but with a recurrent random walk $(S^a_2(t))$ prove
this equality when $q(a)=(1,0)$.

Suppose now that $q(a)\notin\{(1,0), (0,1)\}$. Then by the strong law
of large
numbers, $S^a(t)/t \to m(a)$ almost surely as $t\to\infty$ for any
initial state $S^a(0) = z$. From this, it follows that for
any $S^a(0) = z$ and $\varepsilon> 0$ there is an almost surely finite
positive random variable $N_{z,\varepsilon}$
such that $ |S^a(t) - m(a) t| < \varepsilon t$ for all $t \geq
N_{z,\varepsilon}$.
Since $q(a)\notin\{(1,0), (0,1)\}$, the both coordinates of the mean
vector $m(a)$ are
positive and nonzero and consequently, there exist $N > 0$ and $\hat
\varepsilon> 0$ for which the set
\[
\{ z\in{\mathbb Z}^2 \dvtx |z - m(a) t| < \hat\varepsilon t \mbox{ for
some } t \geq
N \}
\]
is included to the
quadrant ${\mathbb N}^*\times{\mathbb N}^*$. For the initial state
$S^a(0)=0$, from this
it follows that
almost surely $
S^a(t) \in{\mathbb N}^*\times{\mathbb N}^*$ for all $ t \geq\hat{N}
\doteq  \max\{
N_{0,\hat
\varepsilon},
N\}$.
The minimums
\[
\min_{t\in{\mathbb N}} S^a_1(t) \quad\mbox{and}\quad \min_{t\in{\mathbb N}}
S^a_2(t)
\]
are therefore almost surely finite and consequently, for some
$\hat{z}=(\hat{x},\hat{y})\in{\mathbb N}^*\times{\mathbb N}^*$,
\[
{\mathbb P}_{\hat{z}}( \tau^a = +\infty) = {\mathbb P}_0 \Bigl( \min
_{t\in{\mathbb N}} S^a_1(t)
> -\hat
{x} \mbox{ and } \min_{t\in{\mathbb N}}
S^a_2(t) > -\hat{y} \Bigr) > 0.
\]
The last inequality combined with Proposition \ref{Pr3-1} shows that
\[
1 - {\mathbb E}_{\hat{z}}\bigl( \exp\bigl(a\cdot\bigl(S(\tau)-\hat{z}\bigr)\bigr), \tau<
\infty\bigr) =
{\mathbb P}_{\hat{z}}( \tau^a = +\infty) > 0
\]
for some $\hat{z}=(\hat{x},\hat{y})\in{\mathbb N}^*\times{\mathbb
N}^*$. To
complete our
proof, it is now
sufficient to notice that under the hypotheses \hyperlink{hypoH2}{(H2)}, for any $z\in
{\mathbb N}
^*\times{\mathbb N}^*$, the
probability that the random walk $(S^a(t))$ starting at $z$ hits the
point $\hat{z}$
before the first exit from the quadrant ${\mathbb N}^*\times{\mathbb
N}^*$ is nonzero and
consequently, for
some $t = t(z,\hat{z})\in{\mathbb N}$,
\begin{eqnarray*}
&&1 - {\mathbb E}_z\bigl( \exp\bigl(a\cdot\bigl(S(\tau)-z\bigr)\bigr), \tau< \infty\bigr) \\
&&\qquad=
{\mathbb P}_z( \tau
^a =
+\infty) \\
&&\qquad\geq{\mathbb P}_z\bigl(S^a(t) =\hat{z}, \tau^a > t\bigr)
{\mathbb P}_{\hat{z}}( \tau^a = +\infty) > 0.
\end{eqnarray*}
Lemma \ref{Lem4-1} is therefore proved.
\end{pf}
\begin{lemma}\label{Lem4-2} The function
%
%
\begin{equation}\label{e4-1}
z =(x_1,x_2) \to x_2\exp\bigl(a(1,0) \cdot z\bigr) - {\mathbb E}_z\bigl( S_2(\tau
)\exp\bigl(a(1,0)
\cdot S(\tau)\bigr), \tau<
\infty\bigr)\hspace*{-28pt}
\end{equation}
is well defined and nonnegative on ${\mathbb N}^*\times{\mathbb N}^*$.
\end{lemma}
\begin{pf} Indeed, Proposition \ref{Pr3-4} proves that the function
(\ref{e4-1}) is well
defined. To prove that this function is
nonnegative on ${\mathbb N}^*\times{\mathbb N}^*$, let us notice that
by dominated convergence
theorem from Proposition \ref{Pr3-4} it follows that
%
%
\begin{equation}\label{e4-3}
{\mathbb E}_z\bigl( S_2(\tau)\exp\bigl(a\cdot S(\tau)\bigr), \tau<
\infty\bigr)
= \lim_{n\to\infty} {\mathbb E}_z\bigl( S_2(\tau)\exp\bigl(a\cdot S(\tau)\bigr),
\tau\leq n\bigr).\hspace*{-22pt}
\end{equation}
Moreover, the function $z=(x_1,x_2) \to x_2\exp(a(1,0)\cdot z)$ is
harmonic for the random walk $S(t)$ because according to the definition
of the point $a(1,0)$, for any $z=(x_1,x_2)$,
\[
{\mathbb E}_z\bigl(S_2(1) \exp\bigl(a(1,0)\cdot S(1)\bigr)\bigr) -
x_2\exp\bigl(a(1,0)\cdot z\bigr) =
\frac{\partial\varphi(a_1,a_2)}{\partial a_2} \bigg|_{(a_1,a_2)=a(1,0)}
= 0.
\]
Hence, for $a=a(1,0)$, the sequence $S_2(n)\exp(a\cdot S(n))$ is a
martingale relative to the natural filtration of $(S(n))$ and by the
stopping-time theorem, for any $z=(x_1,x_2)\in{\mathbb
N}^*\times{\mathbb N}^*$,
\begin{eqnarray*}
&&{\mathbb E}_z\bigl( S_2(\tau)\exp\bigl(a\cdot S(\tau)\bigr),
\tau\leq n\bigr) \\
&&\qquad= {\mathbb E}_z\bigl( S_2(\tau\wedge n)\exp
\bigl(a\cdot
S(\tau
\wedge
n)\bigr)\bigr) - {\mathbb E}_z\bigl( S_2(n)\exp\bigl(a\cdot S(n)\bigr),
\tau> n\bigr) \\
&&\qquad= x_2\exp(a\cdot z) - {\mathbb E}_z\bigl( S_2(n)\exp
\bigl(a\cdot S(n)\bigr),
\tau> n\bigr) \leq x_2\exp(a\cdot z),
\end{eqnarray*}
where the last relation holds because on the event $\{\tau> n\}$ one
has $ S_2(n) > 0$. The last inequality combined with (\ref{e4-3})
proves that the function (\ref{e4-1}) is nonnegative on ${\mathbb
N}^*\times{\mathbb N}^*$.
\end{pf}
\begin{pf*}{Proof of Proposition \protect\ref{Pr4-1}}
Suppose first that $a\notin \{ a(1,0), a(0,1)\}$. Then by Lemma
\ref{Lem4-1}, the function $ h_a(z) = \exp(a\cdot z) - {\mathbb
E}_z(\exp(a\cdot S(\tau)), \tau< \infty)$ is finite and strictly
positive on ${\mathbb N}^*\times{\mathbb N}^*$. For the homogeneous
random walk $(S(t))$ on ${\mathbb Z}^2$, the exponential function
$z\to\exp (a\cdot z)$ is harmonic and the function
\[
f(z) = {\mathbb E}_z\bigl(\exp\bigl(a\cdot S(\tau)\bigr), \tau< \infty\bigr)
\]
satisfies the equality $
{\mathbb E}_z(f(S(1)) = f(z)$
for all $z\in{\mathbb N}^*\times{\mathbb N}^*$. The function $h_a(z)
= \exp(a\cdot
z) -
f(z)$ satisfies
therefore the equality
\[
{\mathbb E}_z(h_a(S(1))) = h_a(z)
\]
for all $z\in{\mathbb N}^*\times{\mathbb N}^*$. Moreover, for $z\in
{\mathbb Z}\times{\mathbb Z}
\setminus({\mathbb N}
^*\times{\mathbb N}^*)$,
${\mathbb P}_z$-almost surely, $\tau= 0$ and $S(\tau)=z$ from which it
follows that
\[
h_a(z) = \exp(a\cdot z) - {\mathbb E}_z\bigl(\exp\bigl(a\cdot S(\tau)\bigr), \tau<
\infty
\bigr) = 0\qquad
\forall z\in{\mathbb Z}\times{\mathbb Z}\setminus({\mathbb
N}^*\times{\mathbb N}^*).
\]
Since $Z_+(t)$ is killed at the first time $\tau$ when $S(t)$ exits from
${\mathbb N}^*\times{\mathbb N}^*$ and is identical to $S(t)$ for
$t\leq\tau$, we
conclude that the
function $h_a$ is harmonic for the random walk $(Z_+(t))$.
For $a\notin\{a(1,0), a(0,1)\}$, Proposition \ref{Pr4-1} is
therefore proved.

Consider now the case when $a=a(1,0)=(a_1',a_2')$. Then by Lemma \ref
{Lem4-2}, the function $ h_a(z) = x_2\exp(a\cdot z) - {\mathbb E}_z(
S_2(\tau)\exp(a\cdot S(\tau)), \tau< \infty) $ is well defined and
nonnegative on ${\mathbb N}^*\times{\mathbb N}^*$. To prove that this
function is harmonic for the Markov process $(Z_+(t))$ it is sufficient
to notice that
\begin{eqnarray*}
{\mathbb E}_z(h_a(Z_+(1))) & = & {\mathbb E}_z\bigl(h_a(S(1)), \tau> 1\bigr)\\
& = &
{\mathbb E}_z(h_a(S(1))) = h_a(z)\qquad
\forall z\in{\mathbb N}^*\times{\mathbb N}^*,
\end{eqnarray*}
because $h_a(z) = 0$ for all
$z\in{\mathbb Z}\times{\mathbb Z}\setminus({\mathbb N}^*\times
{\mathbb N}^*)$ and $Z_+(1) = S(1)$ whenever
$\tau> 1$.
To prove that the function $h_a$ is strictly positive, we first
notice that
\begin{eqnarray*}
h_a(z) \exp(-a\cdot z) & = & x_2 - {\mathbb E}_z\bigl(
S_2(\tau)\exp\bigl(a\cdot\bigl(S(\tau)-z\bigr)\bigr), \tau= \tau_2
<\infty\bigr)
\\
&&{} - {\mathbb E}_z\bigl(
S_2(\tau)\exp\bigl(a\cdot\bigl(S(\tau)-z\bigr)\bigr), \tau= \tau_1 <
\tau_2 \leq\infty\bigr),
\end{eqnarray*}
where
\[
x_2 - {\mathbb E}_z\bigl( S_2(\tau)\exp\bigl(a\cdot(S(\tau)_z)\bigr), \tau= \tau_2
<\infty\bigr) \geq x_2 > 0,
\]
because on the event $\{\tau= \tau_2\}$ one has $
S_2(\tau) = S_2(\tau_2) \leq0$. Moreover, by Proposition \ref{Pr3-3},
for $a=a(1,0)=(a_1',a_2')$ and any $\delta> 0$ small enough
there is a point $\hat{a} =(\hat{a}_1,\hat{a}_2)\in\partial D$ with
$\hat{a}_1 < a_1'$ and $\hat{a}_2 = a'_2 + \delta$ such that
\begin{eqnarray*}
&&{\mathbb E}_z\bigl( S_2(\tau)\exp\bigl(a\cdot\bigl(S(\tau)-z\bigr)\bigr), \tau= \tau_1 <
\tau_2 \leq\infty\bigr) \\
&&\qquad\leq\frac{1}{\delta} {\mathbb E}_z\bigl( \exp
\bigl(a\cdot
\bigl(S(\tau
)-z\bigr) + \delta S_2(\tau)\bigr), \tau= \tau_1 <
\tau_2 \leq\infty\bigr)
\\
&&\qquad\leq\frac{1}{\delta} \exp\bigl( (\hat{a}-a)\cdot z\bigr) = \frac
{1}{\delta}
\exp\bigl( -(a_1'
-\hat{a}_1) x_1 + \delta x_2\bigr).
\end{eqnarray*}
Since the right-hand side of the last inequality tends to zero as
$x_1\to\infty$, this
proves that $h_a(z) > 0$ for $z=(x_1,x_2)\in{\mathbb N}^*\times
{\mathbb N}^*$ with
$x_2=1$ and $x_1
>0$ large enough. Since by the Harnack
inequality,
\[
h_a(z') \geq h_a(z) {\mathbb P}_{z'}\bigl(Z_+(n) = z \mbox{ for some } n\in
{\mathbb N}\bigr)
\geq0\qquad \forall z,z'\in{\mathbb N}^*\times{\mathbb N}^*,
\]
using \hyperlink{hypoH2}{(H2)}, we conclude that $h_a(z') > 0$ for all $z'\in{\mathbb
N}^*\times
{\mathbb N}^*$.
Proposition \ref{Pr4-1} is therefore proved.
\end{pf*}

\section{Large deviation results}\label{sec5}
In this section, we obtain large deviation results for the family of
scaled processes
and we
deduce from them the logarithmic asymptotics of the Green function. To
get the large
deviation results for scaled processes $(\varepsilon Z_+([t/\varepsilon
]),$ we need
to show that the original nonscaled
process $(Z_+(t))$ satisfies the following communication condition.

\subsection{Communication condition}
\begin{defi}\label{def5-1} A discrete time Markov chain $(\mathcal{Z}(t))$
on a countable state space $E\subset{\mathbb Z}^d$ is
said to satisfy the communication condition on $E_0\subset E$ if there exist
$\theta>0$ and $C>0$ such that for any $z\not=z'$, $z,z'\in E_0$ there
is a sequence of
points $z_0, z_1, \ldots,z_n\in E_0$ with $z_0=z$, $z_n=z'$ and
$n\leq C|z'-z|$ such that
\[
|z_i-z_{i-1}| \leq C \quad\mbox{and}\quad {\mathbb P}_{z_{i-1}}\bigl({\mathcal
Z}(1) = z_i\bigr)
\geq\theta\qquad \forall
i=1,\ldots,n.
\]
\end{defi}
\begin{prop}\label{Pr5-1} Under the hypotheses \textup{\hyperlink{hypoH2}{(H2)}},
the random walk $(Z_+(t))$ satisfies the communication condition on the
hole space ${\mathbb N}^*\times{\mathbb N}^*$.
\end{prop}
\begin{pf}
Indeed, under the hypotheses \hyperlink{hypoH2}{(H2)}, for $\hat{x} =
(1,1)$ and any unit vector $e\in\{(1,0), (0,1)\}$ there is
$n_e\in{\mathbb N}$ such that ${\mathbb P}_{\hat {x}}(Z_+(n_e) = \hat
{x} + e) > 0$. Hence, there are $u_1^e, \ldots,u_{n_e}^e\in{\mathbb
Z}^2$ with $u_1^e+\cdots+ u_{n_e}^e = e$ such that $\mu(u_k^e) > 0$ and
$\hat{x} + u_1^e+\cdots+ u_k^e \in {\mathbb N} ^*\times{\mathbb N}^*$
for all $k\in\{1,\ldots,n_e\}$. Similarly, for any unit vector
$e\in\{(-1,0), (0,-1)\}$ there is $n_e\in{\mathbb N}$ such that $
{\mathbb P}_{\hat{x}- e}(Z_+(n_e) = \hat{x}) > 0$ and consequently,
there are $u_1^e, \ldots,u_{n_e}^e\in{\mathbb Z}^2$ with $u_1^e+\cdots+
u_{n_e}^e = e$ such that $\mu(u_k^e) > 0$ and $\hat{x} - e +
u_1^e+\cdots+ u_k^e \in{\mathbb N} ^*\times{\mathbb N} ^*$ for all
$k\in\{1,\ldots,n_e\}$. This proves that for any $z,z'\in{\mathbb
N}^*\times {\mathbb N} ^*$ with $|z'-z|=1$ there are $n_e\in{\mathbb
N}^*$ and $u_1^e, \ldots,u_{n_e}^e\in {\mathbb Z}^2$ with $z +
u_1^e+\cdots+ u_{n_e}^e = z'$ such that $\mu(u_k^e) > 0$ and $z +
u_1^e+\cdots+ u_k^e \in{\mathbb N}^*\times {\mathbb N}^*$ for all
$k\in\{1,\ldots,n_e\}$ and consequently, the communication condition is
satisfied with
\[
\theta= \min_e \min_{i=1,\ldots,n_e} \mu(u_i^e) > 0
\quad\mbox{and}\quad
C = \max_e \Bigl\{n_e, {\max_{i=1,\ldots,n_e}}
|u_i^e| \Bigr\}.
\]
\upqed\end{pf}

\subsection{Large deviation properties of scaled processes} Before
formulating our large
deviations results, we recall the definition of the sample
path large deviation principle.

Let $D([0,T],{\mathbb R}^2)$ denote the set of all right continuous functions
with left
limits from $[0,T]$ to ${\mathbb R}^2$ endowed with Skorohod metric
(see Billingsley \cite{Billingsley}).
\begin{defi}
(1) A mapping $I_{[0,T]}\dvtx D([0,T],{\mathbb R}^2)\to
[0,+\infty]$ is a good rate function on $D([0,T],{\mathbb R}^2)$ if
for any
$c\geq0$ and
any compact set $V\subset{\mathbb R}^2$, the set
\[
\bigl\{ \phi\in D([0,T],{\mathbb R}^2)\dvtx \phi(0)\in V \mbox{
and } I_{[0,T]}(\phi) \leq c \bigr\}
\]
is compact in $D([0,T],{\mathbb R}^2)$. According to this definition, a good
rate function is lower semi-continuous.

(2) Let $(Z(t))$ be a Markov process on $E\subset{\mathbb Z}^2$ and let
$Z^\varepsilon(t) =\varepsilon Z([t/\varepsilon])$ for $\varepsilon>0$.
When $\varepsilon\to0$, the family of scaled processes
$(Z^\varepsilon(t) =\varepsilon
Z([t/\varepsilon]),
t\in[0,T])$, is said to
satisfy \textit{a sample path large deviation principle} with a rate function
$I_{[0,T]}$ on $D([0,T],{\mathbb R}^2)$ if for any $T>0$ and $z\in
{\mathbb R}^2$
%
%
\begin{equation}\label{e5-1}
\lim_{\delta\to0} \liminf_{\varepsilon\to0} \inf_{z'\in
\varepsilon E \dvtx
|z'-z|<\delta} \varepsilon
\log{\mathbb P}_{[z'/\varepsilon]} \bigl( Z^\varepsilon(\cdot)\in
{\mathcal
O} \bigr) \geq-\inf_{\phi\in{\mathcal O}\dvtx\phi(0)=z}
I_{[0,T]}(\phi)\hspace*{-32pt}
\end{equation}
for every open set ${\mathcal
O}\subset D([0,T],{\mathbb R}^2)$,
and
%
%
\begin{equation}\label{e5-2}
\lim_{\delta\to0} \limsup_{\varepsilon\to0} \sup_{z' \in
\varepsilon E \dvtx
|z'-z|<\delta}
\varepsilon\log{\mathbb P}_{[z'/\varepsilon]} \bigl( Z^\varepsilon(\cdot
)\in
F \bigr) \leq-\inf_{\phi\in F\dvtx\phi(0)=z} I_{[0,T]}(\phi)\hspace*{-32pt}
\end{equation}
for every closed set $F\subset D([0,T],{\mathbb R}^{2})$.
\end{defi}

${\mathbb P}_{[z/\varepsilon]}$ denotes here the distribution of the
Markov process $(Z(t))$ corresponding to the initial state
$Z(0)=[z/\varepsilon]$ where $[z/\varepsilon]$ is the nearest lattice
point to $z/\varepsilon$ in $E$. For $t\in{\mathbb N}$ and
$\varepsilon>0$, we denote by $[t/\varepsilon]$ the integer part of
$t/\varepsilon$.

By Mogulskii's theorem (see \cite{DZ}), under the hypotheses
\hyperlink{hypoH1}{(H1)}--\hyperlink{hypoH3}{(H3)}, the family of
scaled random walks $S^\varepsilon (t)=\varepsilon S([t/\varepsilon])$
satisfies the sample path large deviation principle with a good rate
function
%
%
\begin{equation}\label{e5-3}\qquad
I_{[0,T]}(\phi) = \cases{\displaystyle\int_0^T (\log\varphi)^*(\dot\phi(t))
\,dt, &\quad if
$\phi$ is absolutely continuous,\vspace*{2pt}\cr
+\infty, &\quad otherwise.}
\end{equation}
The convex conjugate $(\log\varphi)^*$ of the function $\log\varphi$
is defined by
\[
(\log\varphi)^*(v) \doteq  \sup_{a\in{\mathbb R}^{2}} \bigl(a\cdot v -
\log
\varphi(a) \bigr).
\]
Under the hypotheses \hyperlink{hypoH4}{(H4)}, $
(\log\varphi)^*(v) = a\cdot v - \log
\varphi(a)$ whenever $v = \nabla\varphi(a)$
because the function $(\log\varphi)$ is convex and differentiable
everywhere in ${\mathbb R}^2$ (see Lemma 2.2.31 of the book of Dembo and
Zeitouni \cite{DZ}).

Consider now the local processes $(Z^1_+(t))$ and $(Z^2_+(t))$. Recall
that $(Z^1_+(t))$
is the random walk on ${\mathbb Z}\times{\mathbb N}^*$ with transition
probabilities
$p_1(z,z') =
\mu(z'-z)$ which is killed at hitting the half-plane ${\mathbb
Z}\times(-{\mathbb N})$.
Similarly, $(Z^2_+(t))$
is the random walk on ${\mathbb N}^*\times{\mathbb Z}$ with transition
probabilities
$p_2(z,z') =
\mu(z'-z)$ which is killed at hitting the half-plane $(-{\mathbb
N})\times{\mathbb Z}$.
The sample path large deviation
principle for
the scaled processes $\varepsilon(Z^1_+([t/\varepsilon]))$ and
$\varepsilon(Z^2_+([t/\varepsilon
]))$ is proved by
Proposition 4.1 of Ignatiouk-Robert \cite{Ignatiouk06}.
\begin{prop}\label{pr5-1} Under the hypotheses
\textup{\hyperlink{hypoH1}{(H1)}--\hyperlink{hypoH3}{(H3)}}, the family
of scaled processes $Z^{\varepsilon,1}_+(t)=\varepsilon
Z^1_+([t/\varepsilon])$ and $Z^{\varepsilon,2}_+(t)=\varepsilon
Z^2_+([t/\varepsilon])$ satisfies the sample path large deviation
principle with the good rate functions
%
%
\begin{equation}\label{e5-4}\qquad
I_{[0,T]}^{1,+}(\phi) = \cases{ \displaystyle\int_0^T (\log\varphi)^*(\dot\phi(t))
\,dt, & if
$\phi$ is absolutely continuous and\cr
& $\phi(t)\in{\mathbb R}\times{\mathbb R}_+$ for all $t\in
[0,T]$,\cr
+\infty, & otherwise,}
\end{equation}
and
%
%
\begin{equation}\label{e5-5}\qquad
I_{[0,T]}^{2,+}(\phi) = \cases{ \displaystyle\int_0^T (\log\varphi)^*(\dot\phi(t))
\,dt, & if
$\phi$ is absolutely continuous and\cr
&$\phi(t)\in{\mathbb R}_+\times{\mathbb R}$ for all $t\in[0,T]$,\cr
+\infty, & otherwise,}
\end{equation}
respectively.
\end{prop}

For the random walk $(Z_+(t))$ killed at the first exit from the
quadrant \mbox{${\mathbb N}^*\times{\mathbb N}^*$}, with the same arguments as
in the proof of Proposition 4.1 of
Ignatiouk-Robert~\cite{Ignatiouk06}, one gets the following statement.
\begin{prop}\label{pr5-2}
Under the hypotheses
\textup{\hyperlink{hypoH1}{(H1)}--\hyperlink{hypoH3}{(H3)}}, the family of
scaled processes $Z^\varepsilon_+(t)=\varepsilon Z_+([t/\varepsilon])$
satisfies the sample path large deviation principle with the good rate
function
%
%
\begin{equation}\label{e5-6}
I_{[0,T]}^+(\phi) = \cases{\displaystyle\int_0^T (\log\varphi)^*(\dot\phi(t))
\,dt, & if
$\phi$ is absolutely continuous and\cr
&$\phi(t)\in{\mathbb R}_+\times{\mathbb R}_+$ for all $t\in
[0,T]$,\cr
+\infty, & otherwise.}\hspace*{-32pt}
\end{equation}
\end{prop}

\subsection{Large deviation estimates of the Green function}

The large deviation properties of scaled processes imply the large
deviation estimates of the Green function. Recall that $G(z,z')$
denotes the Green function of
the homogeneous random walk $(S(t))$, $G_+^i(z,z')$ denotes the Green
function of the Markov
process $(Z_+^i(t))$, for $i=1,2$, and the Green function of the random walk
$(Z_+(t))$ is denoted by $G_+(z,z')$.
\begin{prop}\label{pr5-3} For any $q\in{\mathbb R}^2_+$, $z\in
{\mathbb N}^*\times{\mathbb N}^*$
and any sequences $\varepsilon_n>0$ and
$z_n\in{\mathbb N}^*\times{\mathbb N}^*$ with $\lim_n\varepsilon_n
= 0$ and $\lim_n \varepsilon
_n z_n
= q$, the
following relations hold
%
%
\begin{eqnarray}\label{e5-8}
\lim_{\delta\to0} \liminf_{n\to\infty} \inf_{z\in{\mathbb
N}^*\times{\mathbb N}
^*\dvtx \varepsilon
_n|z| < \delta} \varepsilon_n\log G_+(z,z_n) & \geq & - a(q)\cdot q,
\\
\label{e5-9}
\lim_{\delta\to0} \liminf_{n\to\infty} \inf_{z\in{\mathbb Z}^2\dvtx
\varepsilon
_n|z| <
\delta} \varepsilon_n\log G(z,z_n) & \geq & - a(q)\cdot q
\end{eqnarray}
and for every $i\in\{1,2\}$,
%
%
\begin{equation}\label{e5-10}
\lim_{\delta\to0} \liminf_{n\to\infty} \inf_{z\in{\mathbb
Z}\times{\mathbb N}
^*\dvtx  \varepsilon
_n |z| <
\delta} \varepsilon_n\log G^i_+(z,z_n) \geq- a(q)\cdot q.
\end{equation}
\end{prop}

The proof of this proposition uses the communication condition of
Proposition~\ref{Pr5-1} and the lower large deviation bound
(\ref{e5-1}) for the families of scaled processes
$\varepsilon(Z_+([t/\varepsilon]))$, $\varepsilon(S([t/\varepsilon
]))$, $\varepsilon(Z^1_+([t/\varepsilon ]))$ and
$\varepsilon(Z^2_+([t/\varepsilon]))$, respectively. It is quite
similar to the proof of Proposition 4.2 of Ignatiouk-Robert
\cite{Ignatiouk06}.

\section{Principal part of the renewal equations}\label{sec6}
For $\delta> 0$ and a sequence of points $z_n\in{\mathbb
N}^*\times{\mathbb N}^*$ with $\lim_n|z_n| = +\infty $ and $\lim_n
z_n/|z_n| =q=(q_1,q_2)$, we define the sequence of functions
$\Xi^q_\delta(z,z_n)$ by letting
%
%
\begin{equation}\label{e6-1}
\Xi^q_\delta(z,z_n) = G(z,z_n) - {\mathbb E}_z \bigl(G(S(\tau), z_n), \tau
<\infty, |S(\tau)| <\delta|z_n| \bigr),
\end{equation}
if the coordinates $q_1$ and $q_2$ of the vector $q$ are nonzero. For
$q=(1,0),$ we put
%
%
\begin{equation}\label{e6-2}\qquad
\Xi^q_\delta(z,z_n) = G^1_+(z,z_n) - {\mathbb E}_z \bigl(G^1_+(S(\tau),
z_n), \tau= \tau_1 < \tau_2, |S(\tau )|<\delta |z_n| \bigr)
\end{equation}
and for $q=(0,1),$ we let
%
%
\begin{equation}\label{e6-3}
\Xi^q_\delta(z,z_n) = G^2_+(z,z_n) - {\mathbb E}_z \bigl(G^2_+(S(\tau),
z_n), \tau=\tau_2 < \tau_1, |S(\tau )|<\delta |z_n| \bigr).\hspace*{-28pt}
\end{equation}
Recall that $G(z,z')$ denotes the Green function of the homogeneous
random walk $(S(t))$ on ${\mathbb Z}^2$ having transition probabilities
$p_S(z,z')=\mu(z'-z)$. The Green function of the random walk
$(Z^1_+(t))$ on ${\mathbb Z}\times{\mathbb N}^*$ having a substochastic
transition matrix $(p_1(z,z') = \mu(z'-z), z,z'\in{\mathbb
Z}\times{\mathbb N}^*)$ is denoted by $G^1_+(z,z')$. Similarly,
$G^2_+(z,z')$ denotes the Green function of the random walk
$(Z^2_+(t))$ on ${\mathbb N}^*\times{\mathbb N}^*$ with a
sub-stochastic transition matrix $(p_2(z,z') = \mu(z'-z),
z,z'\in{\mathbb N}^*\times{\mathbb Z})$. The main result of this
section proves that for any $z\in{\mathbb N}^*\times {\mathbb N} ^*$
and $\delta> 0$, the quantity $\Xi^q_\delta(z,z_n)$ represents the
principal part of right-hand side of the renewal equations (\ref{e2-1})
and (\ref{e2-6}) for $z'=z_n$ when ${\lim_n}|z_n| = +\infty$ and $\lim_n
z_n/|z_n| =q$. This is a subject of the following proposition.
\begin{prop}\label{pr6-1}
Under the hypotheses
\textup{\hyperlink{hypoH1}{(H1)}--\hyperlink{hypoH3}{(H3)}}, for any $z\in
{\mathbb N}^*\times{\mathbb N}^*$, $\delta> 0$ and any sequence
$z_n\in{\mathbb N}^*\times{\mathbb N}^*$ with $\lim _n|z_n| = +\infty$
and $\lim_n z_n/|z_n| =q$,
%
%
\begin{equation}\label{e6-4}
\lim_{n\to\infty} G_+(z,z_n)/ \Xi^q_\delta(z,z_n) = 1.
\end{equation}
\end{prop}

To prove this proposition, we need to investigate the function
\[
\lambda_\varepsilon(q,w) = a(w)\cdot w + a(q - w) \cdot(q - w) -
\varepsilon|w|.
\]
\begin{lemma}\label{lem6-1}
Under the hypotheses \textup{\hyperlink{hypoH1}{(H1)}} and
\textup{\hyperlink{hypoH3}{(H3)}}, for any $q\in{\mathcal S}_+^2$ and $\delta
>0$, there is a small $\varepsilon>0$ such that
%
%
\begin{equation}\label{e6-5}
\inf_{w\in{\mathbb R}^2 \dvtx\inf_{\theta>0} |w - \theta q|\geq\delta}
\lambda_\varepsilon
(q,w) > a(q)\cdot q.
\end{equation}
\end{lemma}
\begin{pf}
Under the hypotheses \hyperlink{hypoH1}{(H1)} and
\hyperlink{hypoH3}{(H3)}, the set $D \doteq  \{a\in{\mathbb R}^2 \dvtx
\varphi(a)\leq1\}$ is compact and strictly convex (see
\cite{Hennequin}), and according to the definition (\ref{e1-2}) of the
mapping $q\to a(q)$, the point $a(q)$ is the only point on the boundary
of the set $D$ where the vector $q$ is normal to $D$. For any nonzero
vector $q\in{\mathbb R}^2$, the point $a(q)$ is therefore the only
point in $D$ where the linear function $a\to a\cdot q$ achieves its
maximum over $a\in D$. Hence, for any $w\in{\mathbb R}^2$,
\[
a(w)\cdot w + a(q-w)\cdot(q-w) \geq a(q)\cdot w + a(q)\cdot(q-w) =
a(q)\cdot q,
\]
where the inequality holds with the equality if and only if $ a(w) =
a(q) = a(q-w)$. Since the mapping $w\to a(w)$ from the unit sphere
${\mathcal S}^2$ to $\partial D = \{a\in{\mathbb R}^2 \dvtx
\varphi(a)=1\}$ is one to one, this proves that
%
%
\begin{equation}\label{e6-6}
a(w)\cdot w + a(q-w)\cdot(q-w) > a(q)\cdot q \qquad\mbox{if } w\notin\{
\theta q \dvtx \theta\geq0\}.
\end{equation}
Moreover, the set $D =\{a\in{\mathbb R}^2 \dvtx \varphi(a)\leq1\}$
being compact, the function $w\to a(w)\cdot w + a(q-w)\cdot(q-w)$ is
convex, finite and therefore continuous on ${\mathbb R}^2$. Hence, for
any $R>0$ and $\delta> 0$,
\[
\varepsilon(R,\delta) \doteq  \mathop{\inf_{w\in{\mathbb R}^2 \dvtx
|w|\leq R}}_{\inf _{\theta
>0} |w -
\theta q|\geq\delta} \bigl(a(w)\cdot w + a(q-w)\cdot(q-w) \bigr) - a(q)\cdot q
> 0
\]
and consequently, for $0 < \varepsilon< \varepsilon(R,\delta)/R$,
\[
\inf_{w\in{\mathbb R}^2 \dvtx |w|\leq R, \inf_{\theta>0} |w - \theta
q|\geq \delta} \lambda_\varepsilon(q,w) > a(q)\cdot q.
\]
To get (\ref{e6-5}), it is now sufficient to show that for any
$\varepsilon
>0$ small enough, there is
$R >0$ such that
%
%
\begin{equation}\label{e6-7}
\inf_{w\in{\mathbb R}^2 \dvtx |w|\geq R, \inf_{\theta>0} |w - \theta
q|\geq
\delta}
\lambda_\varepsilon(q,w) > a(q)\cdot q.
\end{equation}
Here, we use the following estimates: for any $w\in{\mathbb R}^2$ and
$q\in{\mathcal S}^2_+$,
%
%
\begin{eqnarray} \label{e6-8}
&&a(w)\cdot w + a(q-w)\cdot(q-w) - a(q)\cdot q \nonumber\\
&&\qquad= \sup_{a\in D} a\cdot w
+ \sup_{a\in
D} a\cdot(q-w) - a(q)\cdot q \nonumber\\[-8pt]\\[-8pt]
&&\qquad\geq a(w)\cdot w +
a(-w)\cdot(q-w) - a(q)\cdot q \nonumber\\
&&\qquad\geq a(w)\cdot w +
a(-w)\cdot(-w) - 2\max_{a\in D} |a|.\nonumber
\end{eqnarray}
The function $\lambda(w) \doteq  a(w)\cdot w +
a(-w)\cdot(-w)$ is continuous and positively homogeneous:
%
%
\begin{equation}\label{e6-9}
\lambda(w) = |w| \lambda(w/|w|).
\end{equation}
Moreover, the same arguments as in the proof of the inequality (\ref
{e6-6}) show that
\[
\lambda(w) > a(w)\cdot w +
a(w)\cdot(-w) = 0 \qquad\mbox{whenever } a(w) \not= a(-w),
\]
and consequently $\lambda(w) > 0$ for all $w\not= 0$. Hence, letting
\[
\varepsilon_0 \doteq \frac{1}{2}\min_{w\in{\mathbb R}^2 \dvtx |w| = 1}
\lambda(w) > 0
\quad\mbox{and}\quad
c \doteq {2\max_{a\in D} }|a|
\]
and using (\ref{e6-9}) at the right-hand side of (\ref{e6-8}), we get
\[
\lambda_\varepsilon(q,w) - a(q)\cdot q \geq2\varepsilon_0 |w| - c -
\varepsilon|w| \geq
\varepsilon_0|w| - \varepsilon|w|
> 0
\]
for all $0<\varepsilon<\varepsilon_0$ and $w\in{\mathbb R}^2$ with
$|w| > c/\varepsilon_0$. The
inequality (\ref{e6-7}) holds
therefore for $R = c/\varepsilon_0$ and $0 < \varepsilon< \varepsilon
_0$, and the
inequality (\ref{e6-5}) is
satisfied for $0 < \varepsilon< \min\{\varepsilon_0, \varepsilon
(c/\varepsilon_0, \delta)\}$.
\end{pf}
\begin{pf*}{Proof of Proposition \protect\ref{pr6-1}}
Let a sequence of points $z_n\in{\mathbb N}^*\times{\mathbb N}^*$ be
such that $\lim
_n|z_n| = +\infty$ and $\lim_n
z_n/|z_n| =q$. Then by Proposition \ref{pr5-3},
\[
\liminf_{n\to\infty} \frac{1}{|z_n|}\log G_+(z,z_n) \geq-
a(q)\cdot
q\qquad \forall z\in{\mathbb N}^*\times{\mathbb N}^*,
\]
and hence, to get (\ref{e6-4}) it is sufficient to show that
%
%
\begin{equation}\label{e6-10}
\limsup_{n\to\infty} \frac{1}{|z_n|}\log\bigl( \Xi^q_\delta(z,z_n) -
G_+(z,z_n) \bigr) < -
a(q)\cdot q\qquad \forall z\in{\mathbb N}^*\times{\mathbb
N}^*.\hspace*{-35pt}
\end{equation}
Moreover, since the quantities $\Xi^q_\delta(z,z_n) - G_+(z,z_n)$ are
decreasing with respect to $\delta>0$, it is sufficient to prove this
relation for small $\delta>0$. For this, the following estimates are
used: for any $\delta
>0$, $z\in{\mathbb N}^*\times{\mathbb N}^*$ and $n\in{\mathbb N}$,
\begin{eqnarray*}
\Xi^q_\delta(z,z_n) - G_+(z,z_n) & = & {\mathbb E}_z \bigl(G(S(\tau), z_n),
\tau
<\infty,
|S(\tau)|\geq\delta|z_n| \bigr) \\
& \leq & \sum_{w\in{\mathbb Z}^2\setminus({\mathbb N}^*\times{\mathbb
N}^*) \dvtx |w|\geq
\delta|z_n|} G(z,w) G(w,z_n),
\end{eqnarray*}
when the coordinates of the vector $q\in{\mathcal S}_+^2$ are nonzero.
Similarly,
\[
\Xi^q_\delta(z,z_n) - G_+(z,z_n) \leq\sum_{w\in(-{\mathbb
N})\times{\mathbb N}^* \dvtx
|w|\geq
\delta|z_n|} G(z,w) G(w,z_n)
\]
for $q=(1,0)$ and
\[
\Xi^q_\delta(z,z_n) - G_+(z,z_n) \leq\sum_{w\in{\mathbb N}^*\times
(-{\mathbb N}) \dvtx
|w|\geq
\delta|z_n|} G(z,w) G(w,z_n)
\]
for $q=(0,1)$. These estimates show that for any $q=(q_1,q_1)\in
{\mathcal
S}_+^2$, $\delta> 0$, $z\in{\mathbb N}^*\times{\mathbb N}^*$ and
$n\in{\mathbb N}$,
%
%
\begin{equation}\label{e6-11}\quad
\Xi^q_\delta(z,z_n) - G_+(z,z_n) \leq\sum_{w\in{\mathbb Z}^2 \dvtx \inf
_{\theta
>0} |w - \theta q|\geq
\kappa\delta|z_n|} G(z,w) G(w,z_n)
\end{equation}
with
\[
\kappa= \cases{1, &\quad if $q=(q_1,q_2)\in\{(1,0), (0,1)\}$, \cr
\min\{q_1, q_2\}, &\quad otherwise.}
\]
Remark
furthermore that for all $a,a'\in\partial D$ and $z, w, z_n\in
{\mathbb Z}^2$
\[
G(z,w)G(w,z_n) = \exp\bigl( - a\cdot(w-z) - a'\cdot(z_n-w)\bigr) G^a(z,w)
G^{a'}(w,z_n),
\]
where $G^a(z,z')$ denotes the Green function of the twisted
random walk $(S^a(t))$ on
${\mathbb Z}^2$ with transition probabilities (\ref{e3-1}). Since
clearly $G^a(z,w) \leq G^a(w,w) = G(0,0)$ and $G^{a'}(w,z_n) \leq
G^{a'}(z_n,z_n) \leq G(0,0)$,
from this it follows that
\begin{eqnarray*}
&&G(z,w)G(w,z_n)\\
&&\qquad\leq\exp\bigl( - a\cdot(w-z) - a'\cdot(z_n-w)\bigr) (G(0,0) )^2\\
&&\qquad\leq\exp\bigl( - a\cdot w- a' \cdot(|z_n|q-w) + a\cdot z - a' \cdot(z_n
- |z_n|q)\bigr) (G(0,0))^2\\
&&\qquad\leq\exp\bigl( - a\cdot w- a' \cdot(|z_n|q-w) + c(|z| + |z_n - |z_n|q |)\bigr)
(G(0,0))^2
\end{eqnarray*}
with $c \doteq  \max_{a\in D} |a|$. Letting moreover $a = a(w/|z_n|) =
a(w)$ and $a' = a(q - w/|z_n|)$ and using the last inequality at the
right-hand side of (\ref{e6-11}), we obtain
\begin{eqnarray*}
&&\bigl(\Xi^q_\delta(z,z_n) - G_+(z,z_n) \bigr)\exp( - c|z| - c |z_n - |z_n|q| )/
(G(0,0) )^2\nonumber
\\
&&\qquad\leq\sum_{w\in{\mathbb Z}^2 \dvtx \inf_{\theta>0} |w - \theta
q|\geq\kappa
\delta|z_n|}
\exp\bigl( - a(w)\cdot w - a(q-w/|z_n|) \cdot(q|z_n| - w) \bigr)
\\
&&\qquad\leq\sum_{w\in{\mathbb Z}^2 \dvtx \inf_{\theta>0} |w - \theta
q|\geq \kappa \delta|z_n| } \exp\bigl( - |z_n|
\lambda_\varepsilon(q,w/|z_n|) - \varepsilon|w| \bigr)
\end{eqnarray*}
for any $\varepsilon> 0$. Since $\lim_n z_n/|z_n|= q$ and the
series $
\sum_{w\in{\mathbb Z}^2} \exp(-\varepsilon|w|)$
converge for every $\varepsilon>0$, from this inequality it follows
that the
left-hand side of
(\ref{e6-10}) does not exceed
\begin{eqnarray*}
&&\limsup_{n\to\infty} \frac{1}{|z_n|}\log\sum
_{w\in{\mathbb Z}
^2 \dvtx\inf_{\theta>0} |w - \theta q|\geq\kappa\delta|z_n|}
\exp \bigl( - |z_n| \lambda_\varepsilon(q, w/|z_n|) - \varepsilon|w| \bigr)
\\
&&\qquad\leq- \inf_{w\in{\mathbb R}^2 \dvtx\inf_{\theta>0} |w - \theta
q|\geq \kappa \delta} \lambda_\varepsilon(q, w).
\end{eqnarray*}
When combined with Lemma \ref{lem6-1}, the last inequality proves
(\ref{e6-10}).
\end{pf*}

\section{Uniform ratio limit theorem for Markov-additive
processes}\label{sec7}
In this section, we improve the ratio limit theorem of the paper \cite
{Ignatiouk06}. This
result is next applied to get the desirable estimates (\ref{e2-5}) and
(\ref{e2-8}) for the
local processes $(S(t))$ and $(Z_+^1(t))$.

\subsection{Uniform ratio limit theorem for general Markov-additive processes}
Recall that a Markov chain ${\mathcal Z}(t)=(A(t),M(t))$ on a countable
set ${\mathbb Z}^{d}\times E$ with transition probabilities
$p((x,y),(x',y'))$ is called \textit{Markov-additive} if
\[
p((x,y),(x',y')) = p\bigl((0,y),(x'-x,y')\bigr) \qquad\mbox{for all $x,x'\in{\mathbb Z}^{d}$,
$y,y'\in E$.}
\]
The first component $A(t)$ of ${\mathcal Z}(t)=(A(t),M(t))$ is said to be
an \textit{additive} part of the process
${\mathcal Z}(t)$, and the second component $M(t)$ is its
\textit{Markovian part}.
The assumptions we need on the Markov-additive process $({\mathcal
Z}(t)=(A(t),M(t)))$ are the
following:

\begin{enumerate}[(A1)]
\item[(A1)]\hypertarget{conditionA1} \mbox{\textit{The Markov chain}} $({\mathcal Z}(t))$ \textit{is irreducible
on} ${\mathbb Z} ^d\times E$.
\item[(A2)]\hypertarget{conditionA2} $E\subset{\mathbb R}^l$ \textit{for
some} $l\in{\mathbb N}$ \textit{and the function}
%
%
\begin{equation}\label{e7-1}
\hat\varphi(a) = \sup_{z\in{\mathbb Z}^{d}\times E}
{\mathbb E}_{z}\bigl( \exp\bigl(a\cdot\bigl({\mathcal Z}(1)-z\bigr)\bigr)\bigr)
\end{equation}
\textit{is finite in a neighborhood of zero in} ${\mathbb R}^{d+l}$.
\end{enumerate}
Remark that the Markov-additive process $({\mathcal Z}(t))$ is not
assumed to be stochastic: its transition matrix can be strictly
substochastic in some points $z = (x,y)\in{\mathbb Z}^{d}\times
E$.

The following property of Markov-additive processes is essential in our
analysis. ${\mathcal G}(z,z')$ denotes here
the Green function of the Markov process $({\mathcal Z}(t))$.
\begin{prop}\label{pr7-1} Let the Markov-additive processes ${\mathcal
Z}(t)=(A(t),M(t))$ be transient and satisfy the conditions
\textup{\hyperlink{conditionA1}{(A1)}} and \textup{\hyperlink{conditionA2}{(A2)}}.
Suppose moreover that for given $w,w'\in{\mathbb Z}^d\times\{0\}$ the
inequality
%
%
\begin{equation}\label{e7-2}
\inf_{z\in{\mathbb Z}^d\times E} \min\bigl\{ {\mathbb P}_z\bigl({\mathcal
Z}(n)=z+w\bigr), {\mathbb P}_z\bigl({\mathcal
Z}(n)=z+w'\bigr) \bigr\} > 0
\end{equation}
holds with some $n>0$. Then for any $0 < \sigma< 1$ and $r > 0$ there
are $C>0$ and
$\theta> 0$ such that
%
%
\begin{eqnarray}\label{e7-3}\quad
{\mathcal G}(z,z') &\leq& \frac{1 + \sigma +C/|z'|}{1 - \sigma} {\mathcal
G}(z + w - w',z')\nonumber\\[-8pt]\\[-8pt]
&&{} + C \exp(- \theta |z'| + r |z|)\nonumber
\end{eqnarray}
for all $z,z'\in{\mathbb Z}^d\times E$.
\end{prop}
\begin{pf} In a particular case, for $n=1$,
this statement was proved in the core of the proof of Proposition 3.2
of the
paper \cite{Ignatiouk06} by using the method of Bernoulli part
decomposition due to Foley
and McDonald \cite{FoleyMcDonald}. When $n > 1$, for the Green function
\[
\tilde{\mathcal G}(z,z') = \sum_{t=0}^\infty{\mathbb P}_z\bigl({\mathcal
Z}(nt) = z'\bigr),\qquad
z,z'\in{\mathbb Z}^d\times E,
\]
of the included Markov chain $\tilde{\mathcal Z}(t) = {\mathcal
Z}(nt)$, this result proves that for any $r > 0$ and $0 < \sigma< 1$
there are $\tilde{C}>0$ and $\tilde\theta> 0$ such that
\begin{eqnarray*}
\tilde{\mathcal G}(z,z') &\leq& \frac{1 + \sigma +\tilde{C}/|z'|}{1 -
\sigma} \tilde{\mathcal G}(z + w - w',z')\\
&&{} + \tilde {C} \exp(-\tilde\theta|z'| + r |z|)\qquad \forall
z,z'\in{\mathbb Z}^d\times E.
\end{eqnarray*}
Since clearly, $\tilde{\mathcal G}(z + w - w',z') = \tilde{\mathcal
G}(z, z' + w' - w)$ for all $w,w'\in{\mathbb Z}^{d}\times\{0\}$ and
\[
{\mathcal G}(z,z') = \sum_{t=0}^{n-1} \sum_{z''\in{\mathbb Z}^d\times
E} {\mathbb P} _z\bigl({\mathcal Z}(t) = z''\bigr) \tilde{\mathcal
G}(z',z'')\qquad
\forall z,z'\in{\mathbb Z}^d\times E,
\]
from this it follows that
\begin{eqnarray*}
{\mathcal G}(z,z') &\leq& \frac{1 + \sigma +\tilde{C}/|z'|}{1 - \sigma}
\sum_{t=0}^{n-1} \sum_{z''\in{\mathbb Z} ^d\times E} {\mathbb
P}_z\bigl({\mathcal Z}(t) = z''\bigr) \tilde{\mathcal G}(z'',z'+w'-w) \\
&&{} +
\tilde{C} \sum_{t=0}^{n-1} \sum_{z''\in{\mathbb Z}^d\times E} {\mathbb
P}_z\bigl({\mathcal Z}(t) = z''\bigr) \exp(- \tilde \theta |z'| + r |z''|)
\\
&\leq& \frac{1 + \sigma+\tilde{C}/|z'|}{1 - \sigma} {\mathcal
G}(z,z'+w'-w)
\\
&&{} + \tilde{C} \exp(-\tilde\theta|z'|) \sum_{t=0}^{n-1}\sum
_{z''\in{\mathbb Z}
^d\times E} {\mathbb P}_z\bigl({\mathcal Z}(t) = z''\bigr) \exp(r
|z''|),
\end{eqnarray*}
where
\begin{eqnarray*}
&&\sum_{z''\in{\mathbb Z}^d\times E} {\mathbb P}_z\bigl({\mathcal Z}(t)
= z''\bigr) \exp(r |z''|)\\
&&\qquad\leq
4 \max_{a\in{\mathbb R}^2 \dvtx
|a|=r} \sum_{z''\in{\mathbb Z}^d\times E} {\mathbb P}_z\bigl({\mathcal
Z}(t) = z''\bigr) \exp
(a\cdot
z'') \\
&&\qquad\leq4 \max_{a\in{\mathbb R}^2 \dvtx |a|\leq r}
\hat\varphi(a)^t\qquad
\forall t\in
{\mathbb N}.
\end{eqnarray*}
When $n>1$, inequality (\ref{e7-3}) holds therefore for $r>0$ small
enough with $\theta= \tilde\theta$ and
\[
C = 4 \tilde{C} \sum_{t=0}^{n-1} \max_{a\in{\mathbb R}^2 \dvtx |a|\leq
r} \hat
\varphi
(a)^t < \infty.
\]
To complete the proof of this proposition, it is now sufficient to
notice that the right-hand side of (\ref{e7-3}) is increasing with respect to $r>0$. Hence,
if the
inequality (\ref{e7-3}) holds with some $C>0$ and $\theta>0$ for a
small $r>0$, then it is
also satisfied for large $r>0$ with the same constants $C$ and $\theta$.
\end{pf}

The following statement is an immediate consequence of Proposition \ref{pr7-1}.
From now on, for the sake of simplicity of expressions, we will use the
following
notation
%
%
\begin{eqnarray}
{\mathop{\underline{\operatorname{Lim}}}_{\delta,n,z}}=\lim_{\delta
\to0} \liminf_{n\to\infty} \inf_{z\in{\mathbb Z}
^{d}\times
E\dvtx |z| < \delta|z_n|},\nonumber\\[-8pt]\\[-8pt]
{\mathop{\overline{\operatorname{Lim}}}_{\delta,n,z}}=\lim_{\delta
\to0} \limsup_{n\to\infty} \sup_{z\in{\mathbb Z}
^{d}\times
E\dvtx |z| < \delta|z_n|}\nonumber.
\end{eqnarray}
\begin{prop}\label{pr7-2} Let a sequence $z_n\in{\mathbb Z}^{d}\times E$
be such that $\lim_n |z_n| = \infty$ and
%
%
\begin{equation}\label{e7-4}
\mathop{\underline{\operatorname{Lim}}}_{\delta,n,z}
\frac{1}{|z_n|} \log{\mathcal G}(z, z_n) \geq0.
\end{equation}
Then under the hypotheses of Proposition \ref{pr7-1},
%
%
\begin{equation}\label{e7-5}
\mathop{\underline{\operatorname{Lim}}}_{\delta,n,z}
\frac{{\mathcal G}(z+w',z_n)}{{\mathcal G}(z+w,z_n)} = {\mathop
{\overline{\operatorname{Lim}}}_{\delta,n,z}}\frac{{\mathcal
G}(z+w',z_n)}{{\mathcal G}(z+w,z_n)} = 1.
\end{equation}
\end{prop}
\begin{pf} Indeed, by Proposition \ref{pr7-1},
for any $r > 0$ and $0 < \sigma< 1$ there are $C>0$ and
$\theta> 0$ such that
\[
{\mathcal G}(z+w',z_n) \leq\frac{1 + \sigma
+C/|z_n|}{1 - \sigma} {\mathcal G}(z + w,z_n) + C \exp(- \theta|z_n|
+ r |z|)
\]
for all $z,z'\in{\mathbb Z}^d\times E$ and consequently,
%
%
\begin{equation}\label{e7-6}
\mathop{\overline{\operatorname{Lim}}}_{\delta,n,z}\frac{{\mathcal
G}(z+w',z_n)}{{\mathcal G}(z+w,z_n)} \leq\frac{1
+ \sigma}{1 - \sigma} + C \,\mathop{\overline{\operatorname
{Lim}}}_{\delta,n,z}\frac{\exp(- \theta|z_n| + r
\delta
|z_n|)}{{\mathcal
G}(z+w,z_n)}.
\end{equation}
Moreover, (\ref{e7-4}) shows that the sequence $\exp(- \theta|z_n| + r
\delta |z_n|)$ tends to zero as $n\to\infty$ faster than the sequence
$1/{\mathcal G}(z+w,z_n)$. From this. it follows that the second term
of the right-hand side of (\ref{e7-6}) is equal to zero and hence,
letting $\sigma\to0$ we conclude that
\[
\mathop{\overline{\operatorname{Lim}}}_{\delta,n,z}{\mathcal
G}(z+w',z_n)/{\mathcal G}(z+w,z_n) \leq1.
\]
To prove the inequality
\[
\mathop{\underline{\operatorname{Lim}}}_{\delta,n,z}{\mathcal
G}(z+w',z_n)/{\mathcal G}(z+w,z_n) \geq1
\]
it is now sufficient to exchange the roles of $w$ and $w'$. The equalities
(\ref{e7-5}) are therefore verified.
\end{pf}

Suppose now that the Markov process $({\mathcal Z}(t))$ satisfies the
communication
condition \ref{def5-1} on ${\mathbb Z}^{d}\times E$.
Then there is a bounded function
$n_0 \dvtx E \to{\mathbb N}^*$ such that for any $z=(x,y)\in{\mathbb
Z}^{d}\times E$,
\[
{\mathbb P}_{(x,y)}\bigl({\mathcal Z}(n_0(y)) = (x,y)\bigr) \geq\theta^{n_0(y)}
> 0
\]
and hence, there is $k \in{\mathbb N}^*$ (for instance, $k = n!$ with
$n =
\max
_y n_0(y)$) such
that
\[
{\mathbb P}_z\bigl({\mathcal Z}(k) = z\bigr) \geq\theta^{k}\qquad \forall z\in
{\mathbb Z}^{d}\times E.
\]
We denote by $\hat{k}$ the greatest common divisor of the set of all
integers $k >0$ for
which
%
%
\begin{equation}\label{e7-7}
\inf_{z\in Z^{d}\times E} {\mathbb P}_z\bigl({\mathcal Z}(k) = z\bigr) > 0.
\end{equation}
The following statement is a refined version of the ratio limit theorem obtained
in~\cite{Ignatiouk06}.
\begin{prop}\label{pr7-3}
Let a Markov-additive process ${\mathcal Z}(t)=(A(t),M(t))$ be
transient and satisfy the communication condition \ref{def5-1} and the
condition \textup{\hyperlink{conditionA2}{(A2)}}. Suppose moreover that
a sequence of points $z_n\in{\mathbb Z}^{d}\times E$ satisfies the
inequality (\ref{e7-4}) with ${\lim_n} |z_n| = \infty$. Then
\[
\mathop{\underline{\operatorname{Lim}}}_{\delta,n,z}
{\mathcal G}(z+\hat{k}w,z_n)/{\mathcal G}(z,z_n) = \mathop{\overline
{\operatorname{Lim}}}_{\delta,n,z}{\mathcal G}(z+\hat
{k}w,z_n)/{\mathcal G}(z,z_n) = 1
\]
for all $w\in{\mathbb Z}^{d}\times\{0\}$.
\end{prop}
\begin{pf}
Indeed, let ${\mathcal K}$ be the set of all integers for which the
inequality (\ref{e7-7}) holds. Because of the communication condition
\ref{def5-1}, for any $w\in{\mathbb Z}\times\{0\}$ there are
$\varepsilon> 0$ and a bounded function $n \dvtx E\to{\mathbb N}^*$
such that
\[
\inf_{z\in{\mathbb Z}^d\times E} {\mathbb P}_z\bigl({\mathcal Z}(n(y)) = z
+ w\bigr) \geq\varepsilon.
\]
Using the Markov property, we get therefore
\[
\inf_{z\in{\mathbb Z}^d\times E} {\mathbb P}_z\bigl({\mathcal Z}(k n(y)) =
z + k w\bigr) \geq\varepsilon^k
\]
for any $k\in{\mathbb N}^*$ and consequently,
\[
\inf_{z\in{\mathbb Z}^d\times E} \min\bigl\{{\mathbb P}_z\bigl({\mathcal Z}(k
n(y)) = z + k w\bigr), {\mathbb P} _z\bigl({\mathcal Z}(k n(y)) = z\bigr) \bigr\} >
0\qquad
\forall k\in{\mathcal K}.
\]
By Proposition \ref{pr7-2}, from this it follows that
%
%
\begin{equation}\label{e7-8}
\mathop{\underline{\operatorname{Lim}}}_{\delta,n,z}
{\mathcal G}(z+ k w,z_n)/{\mathcal G}(z,z_n) =
\mathop{\overline{\operatorname{Lim}}}_{\delta,n,z}{\mathcal G}(z+
k w,z_n)/{\mathcal G}(z,z_n) = 1
\end{equation}
for all $w\in{\mathbb Z}^d\times\{0\}$ and $k\in{\mathcal K}$.
Consider now the
subgroup $\langle{\mathcal K}\rangle$ of ${\mathbb Z}$ generated
by~${\mathcal K}$. Since (\ref{e7-8}) is satisfied for all $w\in
{\mathbb Z}
^d\times\{
0\}$ one can
replace $w$ in the left-hand side of (\ref{e7-8}) by $-w$ and hence,
(\ref{e7-8}) holds
also for any $k\in- {\mathcal K}$. Moreover, if (\ref{e7-8}) is satisfied
for some $k=k_1$
and $k=k_2$ then the same relation is clearly satisfied for
$k=k_1+k_2$. This proves that
(\ref{e7-8}) holds for any $k\in\langle{\mathcal K}\rangle$
and in particular for $k=\hat{k}$ because $\hat{k}\in\langle
{\mathcal
K}\rangle$ (see
Lemma A.1 of Seneta \cite{Seneta}).
\end{pf}

\subsection{Applications to local processes}
According to the above definition, our homogeneous random walk
$(S(t))$ on ${\mathbb Z}^2$ is Markov-additive: its additive part is
the process
$S(t)$ itself
and the Markovian part is empty. The quantity
$\hat{k}$ is here the period of the random walk $(S(t))$. Proposition
\ref{pr7-3} applied
for the process $(S(t))$ with $d=2$ and $E=\varnothing$ and combined with
the estimates
(\ref{e5-9}) yields the following statement.
\begin{prop}\label{pr7-4} For any sequence of points
$z_n\in{\mathbb Z}^2$ with $\lim_n|z_n| = +\infty$ and $\lim_n z_n/|z_n|
=q\in
{\mathcal S}^2$,
%
%
\begin{eqnarray}\label{e7-10}
&&\mathop{\underline{\operatorname{Lim}}}_{\delta,n,z}\exp
\bigl(-a(q)\cdot\hat{k}w\bigr) G(z+\hat{k}w,z_n)/G(z,z_n)
\nonumber\\[-8pt]\\[-8pt]
&&\qquad=\mathop{\overline{\operatorname{Lim}}}_{\delta,n,z}\exp\bigl(-a(q)\cdot
\hat{k}w\bigr) G(z+\hat{k}w,z_n)/G(z,z_n) = 1\nonumber
\end{eqnarray}
for all $w\in{\mathbb Z}^2$.
\end{prop}
\begin{pf}
Indeed, for any $a\in\partial D$, the twisted homogeneous random walk
$(S^a(t))$ defined by (\ref{e3-1}) satisfies the communication
condition \ref{def5-1} and the condition \hyperlink{conditionA2}{(A2)}.
The condition \hyperlink{conditionA2}{(A2)} is satisfied because of the
assumption \hyperlink{hypoH3}{(H3)}, and the communication condition
\ref{def5-1} is satisfied because the random walk $(S^a(t))$ is
irreducible (see the proof of Lemma 4.1 in \cite{Ignatiouk06} for more
details). Moreover, the Green function $G^a(z,z')$ of the twisted
random walk $(S^a(t))$ satisfies the equality
%
%
\begin{equation}\label{e7-11}
G^a(z,z') = G(z,z')\exp\bigl(a\cdot(z'-z)\bigr)\qquad \forall z,z'\in{\mathbb Z}^2.
\end{equation}
Hence, for any sequence of points
$z_n\in{\mathbb Z}^2$ with $\lim_n|z_n| = +\infty$ and $\lim_n z_n/|z_n|
=q\in
{\mathcal S}^2$, using (\ref{e5-9}) we get
\[
\mathop{\underline{\operatorname{Lim}}}_{\delta,n,z}
\frac{1}{|z_n|} \log G^{a(q)}(z, z_n) \geq0
\]
and consequently, by Proposition \ref{pr7-3},
\[
\mathop{\underline{\operatorname{Lim}}}_{\delta,n,z}
G^{a(q)}(z+\hat{k}w,z_n)/G^{a(q)}(z,z_n) = \mathop{\overline
{\operatorname{Lim}}}_{\delta,n,z}G^{a(q)}(z+\hat
{k}w,z_n)/G^{a(q)}(z,z_n) = 1.
\]
The last relations combined with (\ref{e7-11}) prove (\ref{e7-10}).
\end{pf}

We need the following consequence of this proposition.
\begin{cor}\label{cor7-1} Let a sequence of points
$z_n\in{\mathbb Z}^2$ be such that $\lim_n|z_n| = +\infty$ and $\lim
_n z_n/|z_n|
=q\in{\mathcal
S}^2$. Then for any $\sigma>0$ there are $C'>0$, $C''>0$, $\delta>
0$ and $N >0$ such that
%
%
\begin{equation}\label{e7-12}
C' \exp\bigl(a(q)\cdot z - \sigma|z|\bigr) \leq G(z,z_n)/G(0,z_n) \leq C'' \exp
\bigl(a(q)\cdot z + \sigma|z|\bigr)\hspace*{-28pt}
\end{equation}
for all $n \geq N$ and $z\in{\mathbb Z}^2$ with $|z| < \delta|z_n|$.
\end{cor}
\begin{pf} Indeed, the equalities (\ref{e7-10}) show that for any
$\sigma>0$ there are
$\delta> 0$ and $N >0$ such that
%
%
\begin{eqnarray}\label{e7-13}
\exp\bigl(a(q)\cdot\hat{k}e - \hat{k}\sigma/2\bigr) &\leq& {G(u+\hat
{k}e,z_n)}/{G(u,z_n)}\nonumber\\[-8pt]\\[-8pt]
&\leq& \exp\bigl(a(q)\cdot\hat{k}e +
\hat{k}\sigma/2\bigr)\nonumber
\end{eqnarray}
for any unit vector $e\in{\mathbb Z}^2$ and all $n \geq N$, $u\in
{\mathbb Z}^2$ with $|z|
< \delta|z_n|$.
Remark that for any
$z\in{\mathbb Z}^2$ there are unit vectors $e_1\in\{(-1,0), (1,0)\}$ and
$e_2\in\{(0,-1), (0,1)\}$, nonnegative integers $n_1,n_2\in{\mathbb
N}$ and
real numbers
$r_1,r_2 \in[0,1[$ such that
\[
z = \hat{k} (n_1 + r_1) e_1 + \hat{k} (n_2+ r_2) e_2
\quad\mbox{and}\quad
\hat
{k} (n_1 + r_1) + \hat{k} (n_2+ r_2) \leq2|z|.
\]
If $|z| < \delta|z_n|$, then letting $u_0= \hat{k} r_1 e_1 + \hat{k}
r_2 e_2$ and
\[
u_k = \cases{u_0 + k\hat{k}e_1, &\quad for $1\leq k \leq n_1$,\cr
u_0 + n_1\hat{k}e_1 + (k-n_1) \hat{k} e_2, &\quad for $n_1 < k \leq n_1+ n_2$,}
\]
we get $|u_k| \leq(n_1+ r_1) \hat{k} + (n_2 + r_2)\hat{k} \leq2|z| <
2\delta|z_n|$ for all $k=0,\ldots,n_1 +n_2 -1$. The inequalities
(\ref{e7-13}) applied with $u = u_k$ for each $k=0,\ldots , n_1+n_2-1$
prove therefore that
%
%
\begin{eqnarray} \label{e7-14}\hspace*{28pt}
\vienas_{\{|z| < \delta|z_n|\}} \frac{G(z,z_n)}{G(0,z_n)} & \leq &
\frac
{G(u_0,z_n)}{G(0,z_n)}
\prod_{k=0}^{n_1+n_2-1} \vienas_{\{|u_k| < 2\delta|z_n|\}} \frac{G(u_{k+1}
,z_n)}{G(u_{k},z_n)}\nonumber\\
& \leq & \frac{G(u_0,z_n)}{G(0,z_n)} \prod_{k=0}^{n_1+n_2-1}
\exp\bigl(a(q)\cdot(u_k-u_{k-1}) + \hat{k}{\sigma}/{2}\bigr)
\nonumber\\[-8pt]\\[-8pt]
& \leq &
\frac{G(u_0,z_n)}{G(0,z_n)} \exp\bigl(a(q)\cdot(z - z_0) + \hat{k}\sigma
(n_1+n_2)/2\bigr)\nonumber\\
& \leq & \frac{G(u_0,z_n)}{G(0,z_n)} \exp\bigl(a(q)\cdot z + 2\hat{k} |a(q)|
+ \sigma|z|\bigr)\nonumber
\end{eqnarray}
and similarly
%
%
\begin{eqnarray} \label{e7-15}
\vienas_{\{|z| < \delta|z_n|\}}
\frac{G(z,z_n)}{G(0,z_n)}
&\geq&
\vienas_{\{|z| <
\delta|z_n|\}}\frac{G(u_0,z_n)}{G(0,z_n)}\nonumber\\[-8pt]\\[-8pt]
&&{}\times\exp\bigl(a(q)\cdot z - 2\hat{k}|a(q)| -
\sigma|z|\bigr)\nonumber
\end{eqnarray}
for all $n\geq N$.
Remark finally that for any $u\in{\mathbb Z}^2$,
\begin{eqnarray*}
{\mathbb P}_u\bigl(S(t) = 0 \mbox{ for some $t>0$}\bigr) \leq\frac{G(u,z_n)}{G(0,z_n)}
\leq
\frac{1}{{\mathbb P}_0(S(t) = u \mbox{ for some $t>0$})},
\end{eqnarray*}
where ${\mathbb P}_u(S(t) = 0 \mbox{ for some $t>0$}) > 0$ and
${\mathbb P}_0(S(t) = u
\mbox{ for some
$t>0$}) > 0$ because by assumption \hyperlink{hypoH1}{(H1)}, our random
walk $(S(t))$ is irreducible. Using this relation together with (\ref
{e7-14}) and (\ref{e7-15}),
we conclude that (\ref{e7-12}) holds with
\[
C' = \inf_{u\in{\mathbb Z}^2 \dvtx |u|\leq2\hat{k}} {\mathbb
P}_u\bigl(S(t) = 0 \mbox{ for some $t>0$}\bigr) \exp(- 2|a(q)|\hat{k})
\]
and
\[
C'' = \sup_{u\in{\mathbb Z}^2 \dvtx |u|\leq2\hat{k}} \frac{1}{{\mathbb
P}_0(S(t) = u \mbox{ for some $t>0$})} \exp(2|a(q)|\hat{k}).
\]
\upqed\end{pf}

Consider now the random walk $(Z_+^1(t))$ on ${\mathbb Z}\times
{\mathbb N}^*$. Recall that $(Z_+^1(t))$ is identical to $(S(t))$ for $
t < \tau_2 \doteq  \inf\{ n \geq0 \dvtx S(n)\notin{\mathbb Z}\times
{\mathbb N}^*\}$ and killed at the time $\tau_2$. Such a process
$(Z_+^1(t))$ is Markov additive, its additive and Markovian parts are,
respectively, the first and the second coordinates of $(Z_+^1(t))$. To
apply Proposition~\ref {pr7-3} in this case, we need to identify the
greatest common divisor of the set of all integers $k >0$ for which
%
%
\begin{equation}\label{e7-16}
\inf_{z\in Z\times{\mathbb N}^*} {\mathbb P}_z\bigl(Z_+^1(k) = z\bigr) > 0.
\end{equation}
This is a subject of the following lemma.
\begin{lemma}\label{Lem7-1}
The greatest common divisor of the set of all integers $k >0$ for
which (\ref{e7-16}) holds is equal to the period $\hat{k}$ of the
random walk $(S(t))$.
\end{lemma}
\begin{pf}
Indeed, if ${\mathbb P}_0(S(k)=0) > 0$ for some $k\in {\mathbb N}^*$
then there is a sequence of points $u_0,u_1,\ldots, u_k\in{\mathbb
Z}^2$ with $u_0= u_k = 0$ such that
\[
{\mathbb P}_{u_{i-1}}\bigl(S(1) = u_i\bigr) > 0 \qquad\mbox{for all } i=1,\ldots, k.
\]
Moreover, without any restriction of generality one can assume
that for some $l\in\{1,\ldots,
k\}$, the second coordinate of the vectors $e_1=u_1 - u_0,\ldots
,e_l=u_l-u_{l-1}$ is
positive and the second coordinate of the vectors
$e_{l+1}=u_{l+1}-u_l,\ldots, e_k=u_k-u_{k-1}$ is negative or
zero. Then $(0,1) + u_i\in{\mathbb Z}\times{\mathbb N}^*$ for all
$i=0,\ldots, k$ and
\begin{eqnarray*}
\inf_{z\in{\mathbb Z}\times{\mathbb N}^*} {\mathbb P}_z\bigl(Z^{1}_+(k) =
z\bigr) &\geq& {\mathbb P}_{(0,1)}
\bigl(Z^{1}_+(k) = (0,1)\bigr) \\
&\geq&
{\mathbb P}_{(0,1)} \bigl(Z^{1}_+(t)=(0,1) + u_t, \forall t=1,\ldots, k\bigr) \\
&=&
{\mathbb P}_0\bigl(S(t) = u_t, \forall t=1,\ldots, k\bigr) > 0.
\end{eqnarray*}
Since according to the definition of the process $(Z_+^1(t))$,
\[
\inf_{z\in Z\times{\mathbb N}^*} {\mathbb P}_z\bigl(Z_+^1(k) = z\bigr) \leq
{\mathbb P}_z\bigl(S(k) = z\bigr) =
{\mathbb P}_0\bigl(S(k)=0\bigr)\qquad \forall k\in{\mathbb N}^*,
\]
we conclude that (\ref{e7-16}) holds if and only if ${\mathbb
P}_0(S(k)=0) > 0$ and consequently, the greatest common divisor of the
set of all integers $k >0$ for which (\ref{e7-16}) holds is equal to
the period $\hat{k}$ of the random walk $(S(t))$. Lemma \ref{Lem7-1} is
therefore proved.
\end{pf}

From Proposition \ref{pr7-3} applied with $d=1$ and $E={\mathbb N}^*$,
using the estimates (\ref{e5-10}) and Lemma \ref{Lem7-1}, we get the
following statement.
\begin{prop}\label{pr7-5}
For any sequence of points $z_n\in{\mathbb N}^*\times{\mathbb N}^*$
with $\lim_n|z_n| = +\infty $ and $\lim_n z_n/|z_n| =q=(1,0)$,
%
%
\begin{eqnarray}\label{e7-17}
&&\mathop{\underline{\operatorname{Lim}}}_{\delta,n,z}
\exp\bigl(-a(q)\cdot\hat{k}w\bigr) G_+^1(z+\hat{k}w,z_n)/G_+^1(z,z_n)
\nonumber\\[-8pt]\\[-8pt]
&&\qquad=\mathop{\overline{\operatorname{Lim}}}_{\delta,n,z}\exp\bigl(-a(q)\cdot
\hat{k}w\bigr) G_+^1(z+\hat{k}w,z_n)/G_+^1(z,z_n) = 1\nonumber
\end{eqnarray}
for all $w\in{\mathbb Z}\times\{0\}$.
\end{prop}
\begin{pf}
The proof of this proposition is quite similar to the proof of
Proposition \ref{pr7-4}. Proposition \ref{pr7-3} is applied here for
the twisted random walk $(Z_+^{a,1}(t))$ on ${\mathbb Z}\times{\mathbb
N}^*$ with $a=a(q)$, which is identical to $(S^{a}(t))$ for
\[
t < \tau^{a}_2 \doteq  \inf\{ n \geq0 \dvtx S^{a}(n)\notin{\mathbb
Z}\times{\mathbb N} ^*\}
\]
and killed at the time $\tau^a_2$. Lemma 4.1 of \cite{Ignatiouk06}
proves that such a random walk satisfies the communication condition
\ref{def5-1}. The condition \hyperlink{conditionA2}{(A2)} is satisfied
here because by assumption \hyperlink{hypoH3}{(H3)}, for any
$a'\in{\mathbb R}^2$,
\[
\sup_{z\in{\mathbb Z}\times{\mathbb N}^*}
{\mathbb E}_{z}\bigl( \exp\bigl(a'\cdot\bigl(Z_+^{a,1}(1)-z\bigr)\bigr)\bigr) \leq{\mathbb E}_{z}\bigl(
\exp\bigl(a'\cdot
\bigl(S^a(1)-z\bigr)\bigr)\bigr) = \varphi(a'+a) < +\infty.
\]
The greatest common divisor of the set of all integers $k >0$ for
which
\[
\inf_{z\in Z\times{\mathbb N}^*} {\mathbb P}_z\bigl(Z_+^{a,1}(k) = z\bigr) > 0,
\]
is\vspace*{-1pt} clearly the same as for the original process $(Z_+^{1}(t))$. By
Lemma \ref{Lem7-1}, this is the~period $\hat{k}$ of the random walk
$(S(t))$. Finally, the Green function $G^{a,1}_+(z,z')$ of the twisted
random walk $(Z_+^{a,1}(t))$ is related to the Green function
$G^1_+(z,z')$ of the original random walk $(Z_+^{1}(t))$ as follows:
%
%
\begin{equation}\label{e7-18}
G^{a,1}_+(z,z') = G^1_+(z,z')\exp\bigl(a\cdot(z'-z)\bigr)\qquad
\forall z,z'\in {\mathbb Z}^2.
\end{equation}
Using this relation together with (\ref{e5-10}), we conclude that for
any sequence of points
$z_n\in{\mathbb N}^*\times{\mathbb N}^*$ with $\lim_n|z_n| = +\infty
$ and $\lim_n
z_n/|z_n| =q=(1,0)$,
\[
\mathop{\underline{\operatorname{Lim}}}_{\delta,n,z}
\frac{1}{|z_n|} \log G^{a(q),1}_+(z, z_n) \geq0
\]
and consequently, by Proposition \ref{pr7-3}, for any $w\in{\mathbb
Z}\times
\{0\}$,
\begin{eqnarray*}
&&\mathop{\underline{\operatorname{Lim}}}_{\delta,n,z}
G^{a(q),1}_+(z+\hat{k}w,z_n)/G^{a(q),1}_+(z,z_n)\\
&&\qquad = \mathop{\overline{\operatorname{Lim}}}_{\delta,n,z}
G^{a(q),1}_+(z+\hat{k}w,z_n)/G^{a(q),1}_+(z,z_n) =
1.
\end{eqnarray*}
The last relation combined with (\ref{e7-18}) proves (\ref{e7-17}).
\end{pf}

From Proposition \ref{pr7-5}, using the same arguments as in the proof of
Corollary \ref{cor7-1} we get the following statement.
\begin{cor}\label{cor7-2}
Let a sequence of points $z_n\in{\mathbb N}^*\times{\mathbb N}^*$ be
such that $\lim_n|z_n| = +\infty$ and $\lim_n z_n/|z_n| =q=(1,0)$. Then
for any $\sigma>0$ there are $C'>0$, $C'' > 0$, $\delta> 0$ and $N
>0$ such that
\[
C' \exp\bigl(a(q)\cdot w - \sigma|w|\bigr) \leq G_+^1(z+w,z_n)/G_+^1(z,z_n)
\leq C'' \exp\bigl(a(q)\cdot w + \sigma|w|\bigr)
\]
for all $n \geq N$, $z\in{\mathbb Z}\times{\mathbb N}^*$ and $w\in
{\mathbb Z}\times\{0\}$ with $\max\{|z|,|w|\} < \delta|z_n|$.
\end{cor}

For the proof of Theorem \ref{th1}, we need moreover the following
stronger statement.
\begin{prop}\label{pr7-6}
Let a sequence $z_n\in{\mathbb N}^*\times{\mathbb N}^*$ be such that
$\lim_n|z_n| = +\infty$ and $\lim_n z_n/|z_n| =q=(1,0)$. Then for any
$\sigma>0$ there are $C>0$, $\delta> 0$ and $N >0$ such that
\[
G_+^1(z,z_n)/G_+^1(z_0,z_n) \leq C\exp\bigl(a(q)\cdot z + \sigma|z|\bigr)
\]
for all $n \geq N$ and $z\in{\mathbb Z}\times{\mathbb N}^*$ with $|z|
< \delta|z_n|$.
\end{prop}

The proof of this proposition uses Corollary \ref{cor7-2} and the
following results.
\begin{lemma}\label{Lem7-3}
Let $(\xi(t))$ be an irreducible homogeneous random walk on ${\mathbb
Z}$ with a zero mean and a finite variance. Denote $T_{0} \doteq
\inf\{t\geq0 \dvtx \xi(t) \leq0\}$ and let $ T \doteq  \inf\{t\geq0
\dvtx \xi(t) = \xi(0) +1\}$. Then $\lim_{n\to \infty} {\mathbb P}_n(T <
T_0) = 1$.
\end{lemma}
\begin{pf}
Indeed, under the hypotheses of this lemma, $T \doteq  \inf\{ t\geq0
\dvtx \xi(t) = \xi(0) +1\}$ is an almost surely finite stopping time
relative to the natural filtration of $(\xi(t))$ and ${\mathbb
P}_{n+1}(T < T_0) = {\mathbb P}_1 (\xi(t) > - n \mbox{ for all } 0\leq
t\leq T )$ for any $n\in{\mathbb N}$. Hence, by monotone convergence
theorem
\begin{eqnarray*}
\lim_{n\to\infty} {\mathbb P}_{n}(T < T_0) & = & \lim_{n\to\infty}
{\mathbb P}_1
\Bigl(\inf
_{0\leq t \leq T} \xi(t) > - n \Bigr) =
{\mathbb P}_1 \Bigl(\inf_{0\leq t \leq T} \xi(t) >
- \infty\Bigr)\\
& = & \sum_{n=0}^\infty{\mathbb P}_1\Bigl(T = n, \inf_{0\leq t \leq n} \xi
(t) >
- \infty\Bigr) \\
& = & \sum_{n=0}^\infty{\mathbb P}_1(T = n) = 1.
\end{eqnarray*}
\upqed\end{pf}

For the random walk $(S(t))$, this lemma implies the following statement.
\begin{lemma}\label{Lem7-4}
Let $ \hat{\tau} = \inf\{t\geq0 \dvtx S_2(t) = S_2(0)+1\}$. Then for
$a=a(1,0)$, ${\mathbb E}_{(0,k)} ( \exp(a\cdot(S(\hat{\tau }) -
(0,k))), \hat{\tau} < \tau_2 ) \to1$ as $k\to\infty$.
\end{lemma}
\begin{pf}
Indeed, consider the twisted random walk $(S^a(t))$ on ${\mathbb Z} ^2$
having transition probabilities $p_a(z,z') = \exp(a\cdot(z'-z))$ with
$a=a(1,0)$. Then the same arguments as in the proof of Proposition
\ref{Pr3-1} show that
\[
{\mathbb E}_{(0,k)} \bigl( \exp\bigl(a\cdot \bigl(S(\hat{\tau}) -
(0,k)\bigr)\bigr), \hat{\tau} < \tau_2 \bigr) = {\mathbb P}_{(0,k)}
(T^a < T^a_0 )
\]
with $ T^a_0 = \inf\{n\geq0 \dvtx S^{a}_2(t)\leq0\}$ and $ T^a =
\inf\{n\geq0 \dvtx S^{a}_2(t) = S^{a}_2(0)+1\}$. Moreover, for
$a=a(1,0)$, the second coordinate $S^{a}_2(t)$ of $S^{a}(t)$ is a
homogeneous random walk on ${\mathbb Z}$ with zero mean
\[
{\mathbb E}_0(S^{a}_2(1)) = {\mathbb E}_0\bigl(S_2(1)\exp\bigl(a\cdot S(1)\bigr)\bigr) =
\frac{\partial}{\partial a_2}\varphi(a_1, a_2) \bigg|_{(a_1, a_2)
=a(1,0)} = 0
\]
and a finite variance because according to the assumption \hyperlink{hypoH3}{(H3)}, the
jump generating function
\[
\alpha\to{\mathbb E}_0(\exp(\alpha S^{a}_2(1))) = \varphi\bigl(a +
(0,\alpha)\bigr)
\]
of $S^{a}_2(t)$ is finite everywhere in ${\mathbb R}$. Lemma
\ref{Lem7-3} applied with $\xi(t) = S^{a}_2(t)$, $T = T^a$ and $T_0 =
T^a_0$ proves therefore that
\[
\lim_{k\to\infty} {\mathbb E}_{(0,k)} \bigl( \exp\bigl(a\cdot\bigl(S(\hat{\tau}) -
(0,k)\bigr)\bigr), \hat{\tau} < \tau_2 \bigr) = \lim_{k\to\infty} {\mathbb P}_{(0,k)}
(T < T_0 ) = 1.\quad
\]
\upqed\end{pf}
\begin{lemma}\label{Lem7-5}
Under the hypotheses of Lemma \ref{Lem7-4}, for any $\varepsilon> 0$
there are $N_\varepsilon> 0$, $k_\varepsilon>0$ and $\sigma
_\varepsilon> 0$ such that for all $N \geq N_\varepsilon$, $k \geq
k_\varepsilon$ and $0 < \sigma\leq\sigma _\varepsilon$,
\begin{eqnarray*}
&&{\mathbb E}_{(0,k)} \bigl( \exp\bigl(a(1,0)\cdot S(\hat{\tau}) - \sigma
|S_1(\hat\tau)|\bigr),
\hat{\tau} < \tau_2,
|S_1(\hat\tau)| < N \bigr)\\
&&\qquad \geq\exp\bigl(-\varepsilon+ a(1,0)\cdot(0,k)\bigr).
\end{eqnarray*}
\end{lemma}
\begin{pf}
Indeed, for any $x\in{\mathbb Z}$, the sequence $k \to {\mathbb P}
_{(0,k)}(S(\hat{\tau}) = (x,k+1),\break \hat{\tau} < \tau_2)$ is
increasing because for any $k\in {\mathbb N}^*$,
\begin{eqnarray*}
&&{\mathbb P}_{(0,k+1)}\bigl(S(\hat{\tau}) =
(x,k+2), \hat{\tau} < \tau_2\bigr)
\\
&&\qquad\geq {\mathbb P}_{(0,k+1)} \bigl(S(\hat{\tau}) =
(x,k+2) \mbox{ and $S_2(t)>1$ for all $t \leq\hat\tau$} \bigr)
\\
&&\qquad= {\mathbb P}_{(0,k)} \bigl(S(\hat{\tau}) =
(x,k+1) \mbox{ and $S_2(t)>0$ for all $t \leq\hat\tau$} \bigr)
\\
&&\qquad= {\mathbb P}_{(0,k)}\bigl(S(\hat{\tau}) = (x,k+1), \hat{\tau} < \tau_2\bigr).
\end{eqnarray*}
By monotone convergence theorem from this, it follows that
\begin{eqnarray*}
&&{\mathbb E}_{(0,k)} \bigl(\exp\bigl(a(1,0)\cdot\bigl(S(\hat{\tau}) - (0,k)\bigr) -
\sigma
|S_1(\hat\tau)|\bigr), \hat{\tau} < \tau_2,
|S_1(\hat\tau)| < N \bigr) \\
&&\qquad= \sum_{x\in{\mathbb Z}\dvtx |x| < N} \exp
\bigl(a(1,0)\cdot(x,1)
-\sigma|x|\bigr) {\mathbb P}_{(0,k)} \bigl(S(\hat{\tau}) =
(x,k+1), \hat{\tau} < \tau_2 \bigr) \\
&&\qquad\to\sum_{x\in{\mathbb Z}} \exp\bigl(a(1,0)\cdot(x,1)\bigr) \lim_{k\to
\infty}{\mathbb P}_{(0,k)}
\bigl(S(\hat{\tau}) =
(x,k+1), \hat{\tau} < \tau_2 \bigr)
\end{eqnarray*}
as $k\to\infty$, $\sigma\to0$ and ${\mathbb N}\to\infty$.
Moreover, using again
monotone convergence
theorem, we get
\begin{eqnarray*}
&&\sum_{x\in{\mathbb Z}} \exp\bigl(a(1,0)\cdot(x,1)\bigr) \lim_{k\to\infty
}{\mathbb P}_{(0,k)}
\bigl(S(\hat{\tau}) =
(x,k+1), \hat{\tau} < \tau_2 \bigr)
\\
&&\qquad= \lim_{k\to\infty} \sum_{x\in
{\mathbb Z}} \exp
\bigl(a(1,0)\cdot(x,1)\bigr) {\mathbb P}_{(0,k)} \bigl(S(\hat{\tau}) =
(x,k+1), \hat{\tau} < \tau_2 \bigr) \\
&&\qquad= \lim_{k\to\infty} {\mathbb
E}_{(0,k)} \bigl(
\exp\bigl(a(1,0)\cdot\bigl(S(\hat{\tau}) - (0,k)\bigr)\bigr), \hat{\tau} < \tau_2 \bigr).
\end{eqnarray*}
Since by Lemma \ref{Lem7-4}, the right-hand side of the last relation
is equal to $1,$ we
conclude that
\[
{\mathbb E}_{(0,k)} \bigl( \exp\bigl(a(1,0)\cdot\bigl(S(\hat{\tau}) - (0,k)\bigr) -
\sigma
|S_1(\hat
\tau)|\bigr), \hat{\tau} < \tau_2,
|S_1(\hat\tau)| < N \bigr) \to1
\]
as $k\to\infty$, $\sigma\to0$ and ${\mathbb N}\to\infty$, and
consequently, for any $\varepsilon> 0$ there are $N_\varepsilon> 0$,
$k_\varepsilon>0$ and $\sigma _\varepsilon> 0$ such that for all $N
\geq N_\varepsilon$, $k \geq k_\varepsilon$ and $0 < \sigma\leq\sigma
_\varepsilon$,
\[
{\mathbb E}_{(0,k)} \bigl( \exp\bigl(a(1,0)\cdot\bigl(S(\hat{\tau}) - (0,k)\bigr) - \sigma
|S_1(\hat \tau)|\bigr), \hat{\tau} < \tau_2, |S_1(\hat\tau)| < N \bigr)
\geq\exp(-\varepsilon).
\]
Lemma \ref{Lem7-5} is therefore proved.
\end{pf}

Consider now an increasing sequence of stopping times $\hat\tau_k$
defined as follows: $\hat\tau_0 \doteq  0$, $\hat\tau_1 \doteq
\hat\tau$ and $\hat {\tau}_{k} \doteq  \inf\{t\geq\hat\tau_{k-1} \dvtx
S_2(t) = S_2(0)+k\}$ for $k\geq2$. Then from Lemma \ref{Lem7-5}, using
strong Markov property, we obtain the following statement.
\begin{lemma}\label{Lem7-6}
Let $a=a(1,0)$. Then for any $\varepsilon> 0$ there are
$C_\varepsilon>0$, $N_\varepsilon> 0$ and $\sigma _\varepsilon> 0$ such
that for all $N \geq N_\varepsilon$, $0 < \sigma\leq\sigma_\varepsilon$
and $k\geq1$,
%
%
\begin{eqnarray}\label{e7-20}
&&{\mathbb E}_{(0,1)} \bigl( \exp\bigl(a\cdot\bigl(S(\hat{\tau}_k) -(0,1)\bigr) - \sigma
|S_1(\hat\tau
_k)|\bigr), \nonumber\\
&&\hspace*{67.4pt}\hat{\tau}_k < \tau_2,
|S_1(\hat\tau_k)| < N(k-1) \bigr)\\
&&\qquad\geq C_\varepsilon\exp
(-k\varepsilon).\nonumber
\end{eqnarray}
\end{lemma}
\begin{pf} Indeed, by strong Markov property, the left-hand side of the
above inequality is greater than
\begin{eqnarray*}
&&{\mathbb E}_{(0,1)} \Biggl(\prod_{l=1}^{k-1}\exp\bigl(a\cdot\bigl(S\bigl(\hat{\tau
}_{l+1} -
S(\hat{\tau}_l)\bigr) - \sigma\bigl|S_1\bigl(\hat\tau_{l+1} - S_1(\hat\tau_{l})\bigr)\bigr|\bigr),
\hat{\tau}_k < \tau_2,\\
&&\qquad\hspace*{101.1pt}
\bigl|S_1\bigl(\hat\tau_{l+1} - S_1(\hat\tau_{l})\bigr)\bigr| < N, \forall1\leq l \leq k-1
\bigr)\Biggr)\\
&&\qquad \geq\prod_{l=1}^{k-1} {\mathbb E}_{(0,l)} \bigl( \exp\bigl(a\cdot\bigl(S(\hat
{\tau})
- (0,l)\bigr) - \sigma|S_1(\hat\tau)|\bigr), \hat{\tau} < \tau_2,
|S_1(\hat\tau)| < N \bigr)
\end{eqnarray*}
and hence, for any $\varepsilon>0$ with the same quantities
$N_\varepsilon>0$,
$\sigma_\varepsilon>0$ and
$k_\varepsilon>0$ as in Lemma \ref{Lem7-5}, the inequality (\ref{e7-20})
holds for all $N \geq
N_\varepsilon$, $0 < \sigma\leq\sigma_\varepsilon$ and $k\geq1$ with
\begin{eqnarray*}
&&C_\varepsilon= \exp(\sigma_\varepsilon k_\varepsilon) \prod
_{l=1}^{k_\varepsilon-1} {\mathbb E}
_{(0,l)} \bigl(
\exp\bigl(a\cdot\bigl(S(\hat{\tau})
- (0,l)\bigr) - \sigma|S_1(\hat\tau)|\bigr),\\
&&\hspace*{186pt} \hat{\tau} < \tau_2,
|S_1(\hat\tau)| < N \bigr).
\end{eqnarray*}
\upqed\end{pf}
\begin{pf*}{Proof of Proposition \protect\ref{pr7-6}}
Let a sequence of points $z_n\in{\mathbb N}^*\times{\mathbb N}^*$ be
such that ${\lim_n}|z_n| = +\infty$ and $\lim_n z_n/|z_n| =q=(1,0)$. To
simplify the notation, we denote throughout the proof of Proposition
\ref{pr7-6}
\[
a(1,0) = a.
\]
Then by Corollary \ref{cor7-2}, for any $\sigma>0$ there are $C'>0$,
$C''> 0$, $\delta_\sigma> 0$ and $n_\sigma>0$ such that
%
%
\begin{equation}\label{e7-21}\hspace*{25pt}
C'\exp\bigl(a\cdot(x,0) -\sigma|x| \bigr) \leq\frac
{G_+^{1}((x,k),z_n)}{G_+^{1}((0,k),z_n)}
\leq C'' \exp\bigl(a\cdot(x,0) + \sigma|x|\bigr)
\end{equation}
for all those $n \geq n_\sigma$, $x\in{\mathbb Z}$ and $k\in{\mathbb
N}^*$ for which
\[
\max\{|x|,k\} <
\delta_\sigma|z_n|.
\]
Furthermore, recall that the process $(Z^1_+(t))$ is identical to the
homogeneous random walk $(S(t))$ on ${\mathbb Z}^2$ before the first
time when
the second
coordinate $S_2(t)$ of $S(t)$ becomes zero or negative and is killed at
the time $\tau_2 \doteq  \inf\{ t \geq0
\dvtx S_2(t) \leq0\}$. Hence, for any $n\in{\mathbb N}$ and $k\in
{\mathbb N}^*$
\begin{eqnarray*}
\frac{G_+^{1}((0,1),z_n)}{G_+^{1}((0,k),z_n)} & \geq &
\sum_{x\in{\mathbb Z}} {\mathbb P}_{(0,1)} \bigl(S(\hat{\tau}) = (x,k),
\hat{\tau} <
\tau_2 \bigr)
\frac{G_+^{1}((x,k),z_n)}{G_+^{1}((0,k),z_n)} \\
& \geq &\sum_{x\in{\mathbb Z}\dvtx |x|
< N(k-1)} {\mathbb P}_{(0,1)} \bigl(S(\hat{\tau}) = (x,k), \hat{\tau}_k
< \tau
_2 \bigr)
\frac{G_+^{1}((x,k),z_n)}{G_+^{1}((0,k),z_n)}
\end{eqnarray*}
and consequently, for all $n \geq n_\sigma$, $N> 0$ and $k \geq1$
satisfying the
inequalities $0 < N(k-1) < \delta_\sigma|z_n|$ and $1 < k<\delta
_\sigma
|z_n|$, using the first
inequality of (\ref{e7-21}) we get
\begin{eqnarray*}
&&\frac{G_+^{1}((0,1),z_n)}{G_+^{1}((0,k),z_n)} \geq\sum_{x\in
{\mathbb Z}\dvtx |x|
< N(k-1)} C' {\mathbb P}_{(0,1)} \bigl(S(\hat{\tau}) =
(x,k), \hat{\tau} < \tau_2 \bigr)\\
&&\hspace*{143.5pt}{}\times\exp\bigl(a\cdot(x,0)-\sigma|x|\bigr).
\end{eqnarray*}
Moreover, the right-hand side of the above inequality is equal to
\[
C' \exp\bigl(- a\cdot(0,k)\bigr) {\mathbb E}_{(0,1)} \bigl(
\exp\bigl(a\cdot S(\hat {\tau}) - \sigma |S_1(\hat\tau)|\bigr),
\hat{\tau} < \tau_2, |S_1(\hat\tau)| < N(k-1) \bigr)
\]
and hence, using Lemma \ref{Lem7-6}, we conclude that for any
$\varepsilon
>0$, there are $C_\varepsilon>0$,
$N_\varepsilon>0$ and $\sigma_\varepsilon>0$ such that
\[
{G_+^{1}((0,1),z_n)}/{G_+^{1}((0,k),z_n)} \geq C' C_\varepsilon\exp\bigl(
-a\cdot
(0,k-1) - k\varepsilon\bigr),
\]
whenever
\[
0 < \sigma< \sigma_\varepsilon,\qquad n \geq n_\sigma,\qquad 1 \leq
k<\delta_\sigma|z_n| \quad\mbox{and}\quad \delta_\sigma|z_n| > (1-k)N_\varepsilon.
\]
Since $|z_n| \to+\infty$, this proves that for any $\varepsilon>0$
there are
$\hat{C}_\varepsilon>0$,
$\hat\delta_\varepsilon>0$ and $\hat{n}_\varepsilon>0$ such that
%
%
\begin{equation}\label{e7-22}\quad
\vienas_{\{k < \hat\delta_\varepsilon|z_n|\}
}{G_+^{1}((0,k),z_n)}/{G_+^{1}((0,1),z_n)} \leq
\hat{C}_\varepsilon\exp\bigl(a\cdot(0,k) + \varepsilon k\bigr)
\end{equation}
for all $n \geq\hat{n}_\varepsilon$ and $k\in{\mathbb N}^*$.

To complete the proof of our proposition, we combine now the estimates
(\ref{e7-22}) with (\ref{e7-21}). From now on, $\varepsilon>0$ and
$\sigma>0$ are arbitrary and independent from each other. For
$n\geq\max\{n_\sigma,\hat {n}_\varepsilon \}$ and $z=(x,k)\in{\mathbb
Z}\times{\mathbb N}^*$ satisfying the inequalities
$|x|\leq\delta_\sigma|z_n|$ and $k \leq \hat\delta_\varepsilon|z_n|$,
the second inequality of (\ref{e7-21}) together with (\ref{e7-22})
imply that
\begin{eqnarray*}
\frac{G_+^{1}((x,k),z_n)}{G_+^{1}((0,1),z_n)} &\leq&
\frac{G_+^{1}((x,k),z_n)}{G_+^{1}((0,k),z_n)} \times
\frac{G_+^{1}((0,k),z_n)}{G_+^{1}((0,1),z_n)} \\
&\leq&
\hat{C}_\varepsilon C'' \exp
\bigl(a\cdot(x,k) + \sigma|x| + \varepsilon k\bigr)
\end{eqnarray*}
and consequently,
\begin{eqnarray*}
\frac{G_+^{1}((x,k),z_n)}{G_+^{1}(z_0,z_n)} & \leq &
\frac{G_+^{1}(z_0,z_n)}{G_+^{1}((0,1),z_n)} \times\hat
{C}_\varepsilon C''
\exp
\bigl(a\cdot(x,k) +
\sigma|x| + \varepsilon k\bigr)\\
& \leq &\frac{\hat{C}_\varepsilon C''}{{\mathbb P}_{z_0}(Z_+^1(t) = (0,1)
\mbox{ for some $t>0$})} \exp\bigl(a\cdot(x,k) +
\sigma|x| + \varepsilon k\bigr).
\end{eqnarray*}
When $\varepsilon= \delta$, the last inequality proves Proposition \ref
{pr7-6} with $\delta= \min\{\hat\delta_\varepsilon, \delta_\sigma\} > 0
$, $N = \max \{ n_\sigma,\hat{n}_\varepsilon\}
> 0$ and $
C = {C''
\hat{C}_\varepsilon}/{{\mathbb P}_{z_0}(Z_+^1(t) = (0,1) \mbox{ for
some $t>0$})}$.
\end{pf*}

\section{\texorpdfstring{Proof of Theorem \protect\ref{th1}}{Proof of Theorem 1}}\label{sec8}
Let a sequence of point $z_n\in{\mathbb N}^*\times{\mathbb N}^*$ be
such that $\lim
_n|z_n| = +\infty$ and $\lim_n
z_n/|z_n| =q\in{\mathcal S}_+^2$. Recall that by Proposition \ref{pr6-1},
for any
$z\in{\mathbb N}^*\times{\mathbb N}^*$ and $\delta>0$,
\[
\lim_{n\to\infty} G_+(z,z_n)/ \Xi^q_\delta(z,z_n) = 1.
\]
To prove
(\ref{e1-4}) it is therefore sufficient to show that for some $\delta>0$,
%
%
\begin{equation}\label{e8-1}\quad
\lim_{n\to\infty} \Xi^q_\delta(z,z_n)/\Xi^q_\delta(z_0,z_n) =
h_{a(q)}(z)/h_{a(q)}(z_0)\qquad
\forall z\in{\mathbb N}^*\times{\mathbb N}^*.
\end{equation}
Consider first the case when the coordinates $q_1$ and $q_2$ of the
vector $q$ are
nonzero. In this case, the quantities $\Xi^q_\delta(z,z_n)$ are defined
by (\ref{e6-1}),
and to get
(\ref{e8-1}) it is sufficient to show that for some $\delta> 0$,
%
%
\begin{eqnarray}\label{e8-2}\quad
\lim_{n\to\infty} \frac{\Xi^q_\delta(z,z_n)}{G(0,z_n)} & \doteq &
\lim
_{n\to\infty}
\frac{G(z,z_n) }{G(0,z_n)}\nonumber\\
&&{} - \lim_{n\to\infty} {\mathbb E}_z \biggl(\frac
{G(S(\tau),
z_n)}{G(0,z_n)}, \tau<\infty,
|S(\tau)| < \delta|z_n| \biggr) \\
&=& h_{a(q)}(z)\qquad
\forall z\in{\mathbb N}^*\times{\mathbb N}^*.\nonumber
\end{eqnarray}
The proof of this relation uses dominated convergence theorem,
Proposition \ref{Pr3-2}, Corollary \ref{cor7-1} and the
results of Ney and Spitzer \cite{NeySpitzer}. Ney and
Spitzer \cite{NeySpitzer} proved that
\[
\lim_{n\to\infty} G(z,z_n)/G(0,z_n) = \exp\bigl(a(q)\cdot z\bigr)\qquad \forall
z\in{\mathbb Z}^2,
\]
by Corollary \ref{cor7-1},
for any $\sigma>0$, there are $C>0$ and $\delta>0$ for such that
\[
\vienas_{\{|z| <\delta|z_n|\}} G(z,z_n)/G(z_0,z_n) \leq C \exp
\bigl(a(q)\cdot z + \sigma|z|\bigr)
\]
for all $n\in{\mathbb N}$ and $z\in{\mathbb Z}^2\setminus({\mathbb
N}^*\times{\mathbb N}^*)$, and by
Proposition \ref{Pr3-2},
%
%
\begin{equation}\label{e8-3}
{\mathbb E}_z\bigl(\exp\bigl(a(q)\cdot S(\tau) + \sigma|S(\tau)|\bigr) , \tau
<\infty\bigr) <
\infty,
\end{equation}
if $\sigma>0$ is small
enough. Hence, by the dominated convergence theorem,
\[
\lim_{n\to\infty}{\mathbb E}_z \biggl(\frac{G(S(\tau),z_n) }{G(0,z_n)},
\tau
<\infty,
|S(\tau)| <\delta|z_n| \biggr)
= {\mathbb E}_z\bigl(\exp\bigl(a(q)\cdot S(\tau)\bigr), \tau<\infty\bigr)
\]
and consequently, (\ref{e8-2}) holds.
When the coordinates of $\lim_n
z_n/|z_n| =q$ are nonzero, the equality (\ref{e8-1}) is therefore
proved.

Suppose now that $\lim_n z_n/|z_n| =q = (1,0)$. For such a vector $q$,
the quantities
$\Xi^q_\delta(z,z_n)$ are defined by (\ref{e6-2}), and to get (\ref
{e8-1}) it is sufficient
to show that for some $\delta>0$ and $C_0 > 0$,
%
%
\begin{eqnarray}\label{e8-5}
&&\lim_{n\to\infty}
\frac{G_+^1(z,z_n) }{G_+^1(z_0,z_n)} - \lim_{n\to\infty}{\mathbb E}_z
\biggl(\frac
{G_+^1(S(\tau),z_n)
}{G_+^1(z_0,z_n)}, \tau= \tau_1 < \tau_2, |S(\tau)| < \delta|z_n|
\biggr)\hspace*{-25pt}
\nonumber\\[-8pt]\\[-8pt]
&&\qquad= C_0 h_{a(q)}(z)\qquad
\forall z\in{\mathbb N}^*\times{\mathbb N}^*.\hspace*{-25pt}\nonumber
\end{eqnarray}
The proof of this equality uses the same arguments as above but with
the help of
Propositions \ref{Pr3-3}, \ref{pr7-6} and the
results of \cite{Ignatiouk06}. Theorem 1 of \cite
{Ignatiouk06} proves the
point-wise convergence
\[
\lim_{n\to\infty} G^1_+(z,z_n)/G^1_+(z_0,z_n) =
h_{a(q),+}^1(z)/h_{a(q),+}^1(z_0)
\]
with a strictly positive function $h_{a(q),+}^1$ on ${\mathbb Z}\times
{\mathbb N}^*$
defined by
%
%
\begin{equation}\label{e8-6}\qquad
h_{a(q),+}^1(z) = x_2\exp\bigl(a(q)\cdot z\bigr) - {\mathbb E}_z
\bigl(S_2(\tau_2)\exp\bigl(a(q)\cdot S(\tau_2)\bigr), \tau_2
<\infty\bigr).
\end{equation}
By Proposition \ref{Pr3-3}, (\ref{e8-3}) holds
if $\sigma>0$ is small enough and by Proposition \ref{pr7-6}, for any
$\sigma>0$ there are $C>0$ and
$\delta>0$ such that
\[
\vienas_{\{|z| < \delta|z_n|\}} G^1_+(z,z_n)/G^1_+(z_0,z_n) \leq C
\exp \bigl(a(q)\cdot z + \sigma|z|\bigr)
\]
for all $n\in{\mathbb N}$ and $z\in{\mathbb Z}\times{\mathbb N}^*$.
By dominated convergence
theorem, from this it
follows that the left-hand
side of (\ref{e8-5}) is equal to
\[
\frac{1}{h_{a(q),+}^1(z_0)} \bigl(h_{a(q),+}^1(z) -
{\mathbb E}_z \bigl(h_{a(q),+}^1(S(\tau)), \tau= \tau_1 < \tau_2 \bigr) \bigr).
\]
Finally, for any $z=(x_1,x_2)\in{\mathbb N}^*\times{\mathbb N}^*$,
from (\ref{e8-6}) it
follows that
\begin{eqnarray*}
&&h_{a(q),+}^1(z) - {\mathbb E}_z\bigl( h_{a(q),+}^1(S(\tau)), \tau=\tau_1 <
\tau_2 \bigr) \\
&&\qquad = x_2\exp\bigl(a(q)\cdot z\bigr) -
{\mathbb E}_z \bigl(S_2(\tau_2)\exp\bigl(a(q)\cdot S(\tau_2)\bigr), \tau_2 <
\infty\bigr) \\
&&\qquad\quad{} - {\mathbb E}_z\bigl( S_2(\tau) \exp\bigl(a(q)\cdot S(\tau)\bigr), \tau= \tau_1
< \tau
_2 \bigr)\\
&&\qquad\quad{} + {\mathbb E}_z \bigl({\mathbb E}_{S(\tau)} \bigl(S_2(\tau_2)\exp\bigl(a(q)\cdot
S(\tau_2)\bigr), \tau_2
<\infty\bigr), \tau= \tau_1 < \tau_2 \bigr).
\end{eqnarray*}
By strong Markov property, the last term of the right-hand side of this
relation is equal to
\begin{eqnarray*}
&&{\mathbb E}_z \bigl(S_2(\tau_2)\exp\bigl(a(q) \cdot
S(\tau_2)\bigr), \tau_1 < \tau_2 < \infty\bigr) \\
&&\qquad= {\mathbb E}_z \bigl(S_2(\tau
_2)\exp\bigl(a(q)
\cdot
S(\tau_2)\bigr), \tau_2 < \infty\bigr)\\
&&\qquad\quad\hspace*{0pt}{} - {\mathbb E}_z \bigl(S_2(\tau_2)\exp
\bigl(a(q)\cdot
S(\tau_2)\bigr), \tau= \tau_2 \leq\tau_1 \bigr)
\end{eqnarray*}
from which if follows that
\begin{eqnarray*}
&&h_{a(q),+}^1(z) -{\mathbb E}_z\bigl( h_{a(q),+}^1(S(\tau)), \tau=\tau_1 <
\tau_2 \bigr) \\
&&\qquad = x_2\exp\bigl(a(q)\cdot z\bigr) -
{\mathbb E}_z \bigl(S_2(\tau)\exp\bigl(a(q)\cdot S(\tau)\bigr), \tau<\infty\bigr) \\
&&\qquad
= h_{a(q)}(z)
\end{eqnarray*}
and consequently, the left-hand side of (\ref{e8-5}) is equal to
$h_{a(q)}(z)/h_{a(q),+}^1(z_0)$.
The equality (\ref{e8-5}) holds therefore with $C_0 =
1/h_{a(q),+}^1(z_0) > 0$ and
hence, for $q=(1,0)$, the equality (\ref{e8-1}) is also proved.

The proof of (\ref{e8-1}) for $q= (0,1)$ uses exactly the same
arguments as above, it is sufficient to exchange the roles
of the first and the second coordinates.

%

%
\printaddresses


\begin{thebibliography}{20}

\bibitem{AliliDoney}
%
\begin{barticle}[mr]
\bauthor{\bsnm{Alili},~\bfnm{L.}\binits{L.}} \AND
\bauthor{\bsnm{Doney},~\bfnm{R.~A.}\binits{R.~A.}}
(\byear{2001}).
\btitle{Martin boundaries associated with a killed random walk}.
\bjournal{Ann. Inst. H. Poincar\'e Probab. Statist.}
\bvolume{37}
\bpages{313--338}.
\bid{doi={10.1016/S0246-0203(00)01069-4}, mr={1831986}}
\end{barticle}
%
\endbibitem

\bibitem{Billingsley}
%
\begin{bbook}[mr]
\bauthor{\bsnm{Billingsley},~\bfnm{Patrick}\binits{P.}}
(\byear{1968}).
\btitle{Convergence of Probability Measures}.
\bpublisher{Wiley}, \baddress{New York}.
\bid{mr={0233396}}
\end{bbook}
%
\endbibitem

\bibitem{Cartier}
%
\begin{bincollection}[mr]
\bauthor{\bsnm{Cartier},~\bfnm{P.}\binits{P.}}
(\byear{1972}).
\btitle{Fonctions harmoniques sur un arbre}.
In \bbooktitle{Symposia {M}athematica, {V}ol. {IX} ({C}onvegno di {C}alcolo
delle {P}robabilit\`a, {INDAM}, {R}ome, 1971)}
\bpages{203--270}.
\bpublisher{Academic Press}, \baddress{London}.
\bid{mr={0353467}}
\end{bincollection}
%
\endbibitem

\bibitem{Chow}
%
\begin{barticle}[mr]
\bauthor{\bsnm{Chow},~\bfnm{Y.~S.}\binits{Y.~S.}}
(\byear{1986}).
\btitle{On moments of ladder height variables}.
\bjournal{Adv. in Appl. Math.}
\bvolume{7}
\bpages{46--54}.
\bid{doi={10.1016/0196-8858(86)90005-9}, mr={834219}}
\end{barticle}
%
\endbibitem

\bibitem{DZ}
%
\begin{bbook}[mr]
\bauthor{\bsnm{Dembo},~\bfnm{Amir}\binits{A.}} \AND
\bauthor{\bsnm{Zeitouni},~\bfnm{Ofer}\binits{O.}}
(\byear{1998}).
\btitle{Large Deviations Techniques and Applications},
\bedition{2nd} ed.
\bseries{Applications of Mathematics (New York)}
\bvolume{38}.
\bpublisher{Springer}, \baddress{New York}.
\bid{mr={1619036}}
\end{bbook}
%
\endbibitem

\bibitem{FoleyMcDonald}
%
\begin{barticle}[mr]
\bauthor{\bsnm{Foley},~\bfnm{Robert~D.}\binits{R.~D.}} \AND
\bauthor{\bsnm{McDonald},~\bfnm{David~R.}\binits{D.~R.}}
(\byear{2005}).
\btitle{Bridges and networks: Exact asymptotics}.
\bjournal{Ann. Appl. Probab.}
\bvolume{15}
\bpages{542--586}.
\bid{doi={10.1214/105051604000000675}, mr={2114982}}
\end{barticle}
%
\endbibitem

\bibitem{Doney02}
%
\begin{barticle}[mr]
\bauthor{\bsnm{Doney},~\bfnm{R.~A.}\binits{R.~A.}}
(\byear{1998}).
\btitle{The {M}artin boundary and ratio limit theorems for killed random
walks}.
\bjournal{J. Lond. Math. Soc. (2)}
\bvolume{58}
\bpages{761--768}.
\bid{doi={10.1112/S0024610798006826}, mr={1678162}}
\end{barticle}
%
\endbibitem

\bibitem{Doob}
%
\begin{barticle}[mr]
\bauthor{\bsnm{Doob},~\bfnm{J.~L.}\binits{J.~L.}}
(\byear{1959}).
\btitle{Discrete potential theory and boundaries}.
\bjournal{J. Math. Mech.}
\bvolume{8}
\bpages{433--458}.
\bid{mr={0107098}}
\end{barticle}
%
\endbibitem

\bibitem{Hennequin}
%
\begin{barticle}[mr]
\bauthor{\bsnm{Hennequin},~\bfnm{Paul-Louis}\binits{P.-L.}}
(\byear{1963}).
\btitle{Processus de {M}arkoff en cascade}.
\bjournal{Ann. Inst. H. Poincar\'e}
\bvolume{18}
\bpages{109--195}.
\bid{mr={0164373}}
\end{barticle}
%
\endbibitem

\bibitem{Hunt}
%
\begin{barticle}[mr]
\bauthor{\bsnm{Hunt},~\bfnm{G.~A.}\binits{G.~A.}}
(\byear{1960}).
\btitle{Markoff chains and {M}artin boundaries}.
\bjournal{Illinois J. Math.}
\bvolume{4}
\bpages{313--340}.
\bid{mr={0123364}}
\end{barticle}
%
\endbibitem

\bibitem{Ignatiouk08}
%
\begin{barticle}[vtex]
\bauthor{\bsnm{Ignatiouk-Robert},~\bfnm{Irina}\binits{I.}}
(\byear{2009}).
\btitle{Martin boundary of a reflected random walk on a
  half-space}.
\bjournal{Probab. Theory Related Fields}.
DOI: \href{http://dx.doi.org/10.1007/s00440-009-0228-4}{10.1007/s00440-009-0228-4}.
\end{barticle}
%
\endbibitem

\bibitem{Ignatiouk06}
%
\begin{barticle}[mr]
\bauthor{\bsnm{Ignatiouk-Robert},~\bfnm{Irina}\binits{I.}}
(\byear{2008}).
\btitle{Martin boundary of a killed random walk on a half-space}.
\bjournal{J.~Theoret. Probab.}
\bvolume{21}
\bpages{35--68}.
\bid{doi={10.1007/s10959-007-0100-3}, mr={2384472}}
\end{barticle}
%
\endbibitem

\bibitem{IMS}
%
\begin{barticle}[mr]
\bauthor{\bsnm{Ignatyuk},~\bfnm{I.~A.}\binits{I.~A.}},
\bauthor{\bsnm{Malyshev},~\bfnm{V.~A.}\binits{V.~A.}} \AND
\bauthor{\bsnm{Shcherbakov},~\bfnm{V.~V.}\binits{V.~V.}}
(\byear{1994}).
\btitle{The influence of boundaries in problems on large deviations}.
\bjournal{Uspekhi Mat. Nauk}
\bvolume{49}
\bpages{43--102}.
\bid{doi={10.1070/RM1994v049n02ABEH002204}, mr={1283135}}
\end{barticle}
%
\endbibitem

\bibitem{KurkovaMalyshev}
%
\begin{barticle}[mr]
\bauthor{\bsnm{Kurkova},~\bfnm{I.~A.}\binits{I.~A.}} \AND
\bauthor{\bsnm{Malyshev},~\bfnm{V.~A.}\binits{V.~A.}}
(\byear{1998}).
\btitle{Martin boundary and elliptic curves}.
\bjournal{Markov Process. Related Fields}
\bvolume{4}
\bpages{203--272}.
\bid{mr={1641546}}
\end{barticle}
%
\endbibitem

\bibitem{Martin}
%
\begin{barticle}[mr]
\bauthor{\bsnm{Martin},~\bfnm{Robert~S.}\binits{R.~S.}}
(\byear{1941}).
\btitle{Minimal positive harmonic functions}.
\bjournal{Trans. Amer. Math. Soc.}
\bvolume{49}
\bpages{137--172}.
\bid{mr={0003919}}
\end{barticle}
%
\endbibitem

\bibitem{NeySpitzer}
%
\begin{barticle}[mr]
\bauthor{\bsnm{Ney},~\bfnm{P.}\binits{P.}} \AND
\bauthor{\bsnm{Spitzer},~\bfnm{F.}\binits{F.}}
(\byear{1966}).
\btitle{The {M}artin boundary for random walk}.
\bjournal{Trans. Amer. Math. Soc.}
\bvolume{121}
\bpages{116--132}.
\bid{mr={0195151}}
\end{barticle}
%
\endbibitem

\bibitem{KilianRaschel}
%
\begin{bmisc}[auto:SpringerTagBib|2009-01-14|16:51:27]
\bauthor{\bsnm{Raschel},~\bfnm{K.}\binits{K.}}
(\byear{2009}).
\bhowpublished{Random walks in the quarter plane absorbed at
the boundary: Exact and asymptotic. Preprint. Available at}
\url{http://arxiv.org/abs/math.PR/0902.2785}.
\end{bmisc}
%
\endbibitem

\bibitem{Rogers}
%
\begin{bbook}[mr]
\bauthor{\bsnm{Rogers},~\bfnm{L.~C.~G.}\binits{L.~C.~G.}} \AND
\bauthor{\bsnm{Williams},~\bfnm{David}\binits{D.}}
(\byear{1994}).
\btitle{Diffusions, {M}arkov Processes, and Martingales, {V}ol. 1.
Foundations},
\bedition{2nd} ed.
\bpublisher{Wiley}, \baddress{Chichester}.
\bid{mr={1331599}}
\end{bbook}
%
\endbibitem

\bibitem{Seneta}
%
\begin{bbook}[mr]
\bauthor{\bsnm{Seneta},~\bfnm{E.}\binits{E.}}
(\byear{1981}).
\btitle{Nonnegative Matrices and {M}arkov Chains},
\bedition{2nd} ed.
\bpublisher{Springer}, \baddress{New York}.
\bid{mr={719544}}
\end{bbook}
%
\endbibitem

\bibitem{Woess}
%
\begin{bbook}[mr]
\bauthor{\bsnm{Woess},~\bfnm{Wolfgang}\binits{W.}}
(\byear{2000}).
\btitle{Random Walks on Infinite Graphs and Groups}.
\bseries{Cambridge Tracts in Mathematics}
\bvolume{138}.
\bpublisher{Cambridge Univ. Press}, \baddress{Cambridge}.
\bid{mr={1743100}}
\end{bbook}
%
\endbibitem

\end{thebibliography}
\end{document}